\documentclass[10pt]{article}%
\usepackage{amsmath}
\usepackage{amsfonts}
\usepackage{mathrsfs}
\usepackage{amssymb, color}
\usepackage[linkcolor=black,anchorcolor=black,citecolor=black]{hyperref}
\usepackage{graphicx}
\numberwithin{equation}{section}
\usepackage[body={15.5cm,21cm}, top=3cm]{geometry}%
\providecommand{\U}[1]{\protect\rule{.1in}{.1in}}
\providecommand{\U}[1]{\protect \rule{.1in}{.1in}}
\newtheorem{theorem}{Theorem}[section]

\newtheorem{definition}[theorem]{Definition}
\newtheorem{example}[theorem]{Example}

\newtheorem{lemma}[theorem]{Lemma}

\newtheorem{proposition}[theorem]{Proposition}
\newtheorem{remark}[theorem]{Remark}

\newenvironment{proof}[1][Proof]{\noindent \textbf{#1.} }{\  \rule{0.5em}{0.5em}}
\DeclareMathOperator*{\esssup}{ess\,sup}
\DeclareMathOperator*{\essinf}{ess\,inf}
\begin{document}
\title{Reflected Solutions of BSDEs Driven by $G$-Brownian Motion}
\author{ Hanwu Li\thanks{School of Mathematics, Shandong University,
lihanwu@mail.sdu.edu.cn.}
\and Shige Peng\thanks{School of Mathematics and Qilu Institute of Finance, Shandong University,
peng@sdu.edu.cn. Li and Peng's research was
partially supported by NSF (No. 10921101) and by the 111
Project (No. B12023).}}
\maketitle
\begin{abstract}
In this paper, we study the reflected solutions of one-dimensional  backward stochastic differential equations driven by $G$-Brownian motion (RGBSDE for short). The reflection keeps the solution above a given stochastic process. In order to derive the uniqueness of reflected $G$-BSDEs, we apply a ``martingale condition" instead of the Skorohod condition. Similar to the classical case, we prove the existence by approximation via penalization.
\end{abstract}


\textbf{MSC-classification}: 60H10, 60H30
\section{Introduction}
El Karoui, Kapoudjian, Pardoux, Peng and Quenez \cite{KKPPQ} studied the problem of BSDE with reflection, which means that the solution to a BSDE is required to be above a certain given continuous boundary process, called the obstacle. For this purpose, an additional continuous increasing process should be added in the equation. Furthermore, this additional process should be chosen in a minimal way so that the Skorohod condition is satisfied. An important observation is that the solution is the value function of an optimal stopping problem.

Due to the importance in BSDE theory and in applications, the reflected problem has attracted a great deal of attention since 1997. Many scholars tried to relax the conditions on the generator and the obstacle process.  Hamadene \cite{H} and Lepeltier and Xu \cite{LX} gave a generalized Skorohod condition and proved the existence and uniqueness when the obstacle process is no longer continuous. Cvitanic and Karaztas \cite{CK} and Hamadene and Lepeltier \cite{HL} studied the case of two reflecting obstacles. They also established the connection between this problem and Dynkin games. Matoussi \cite{M} and Kobylanski, Lepeltier, Quenez and Torres \cite{KLQT} extended the results to the case where the generator is not a Lipschitz function.

We should point out that the classical BSDEs can only provide probabilistic interpretation for quasilinear partial differential equations (PDE for short). Besides, this BSDE cannot be applied to price path-dependent contingent claims in the uncertain volatility model (UVM for short). Motivated by these facts, Peng \cite{P04,P05} systemetically introduced a time-consistent fully nonlinear expectation theory. One of the most important cases is the $G$-expectation theory (see \cite{P10} and the reference therein). In this framework, a new type of Brownian motion and the corresponding stochastic calculus of It\^{o}'s type were constructed. It has been widely used to study the problems of model uncertainty, nonlinear stochastic dynamical systems and fully nonlinear PDEs.

The backward stochastic differential equations driven by $G$-Brownian motion (i.e., $G$-BSDE) can be written in the following way
\begin{displaymath}
Y_t=\xi+\int_t^Tf(s,Y_s,Z_s)ds+\int_t^Tg(s,Y_s,Z_s)d\langle B\rangle_s-\int_t^T Z_s dB_s-(K_T-K_t).
\end{displaymath}
The solution of this equation consists of a triplet of processes $(Y,Z,K)$. The existence and uniqueness of the solution are proved in \cite{HJPS1}. In \cite{HJPS2} the comparison theorem, Feymann-Kac formula and some related topics associated with this kind of $G$-BSDEs were established.

In this paper, we study the case where the solution of a $G$-BSDE is required to stay above a given stochastic process, called the lower obstacle. An increasing process should be added in this equation to push the solution upwards, so that it may remain above the obstacle. According to the classical case studied by \cite{KKPPQ}, we may expect that the solution of reflected $G$-BSDE is a quadruple $\{(Y_t,Z_t,K_t,A_t),0\leq t\leq T\}$ satisfying

\begin{description}
\item[(1)] $Y_t=\xi+\int_t^T f(s, Y_s,Z_s) ds+\int_t^Tg(s,Y_s,Z_s)d\langle B\rangle_s-\int_t^T Z_s dB_s-(K_T-K_t)+A_T-A_t$;
\item[(2)] $(Y,Z,K)\in\mathfrak{S}_G^{\alpha}(0,T)$ and $Y_t\geq S_t$, $0\leq t\leq T$;
\item[(3)] $\{A_t\}$ is continuous and increasing, $A_0=0$ and $\int_0^T (Y_t-S_t)dA_t=0$.
\end{description}

The shortcoming is that the solution of the above problem is not unique.  Thus, to get the uniqueness of the reflected $G$-BSDE, we should reformulate this problem as the following. A triplet of processes $(Y,Z,A)$ is called a solution of reflected $G$-BSDE if the following properties hold: 
\begin{description}
\item[(a)]$(Y,Z,A)\in\mathcal{S}_G^{\alpha}(0,T)$ and $Y_t\geq S_t$;
\item[(b)]$Y_t=\xi+\int_t^T f(s,Y_s,Z_s)ds+\int_t^T g(s,Y_s,Z_s)d\langle B\rangle_s
-\int_t^T Z_s dB_s+(A_T-A_t)$;
\item[(c)]$\{-\int_0^t (Y_s-S_s)dA_s\}_{t\in[0,T]}$ is a non-increasing $G$-martingale.
\end{description}
Here, we denote by $\mathcal{S}_G^{\alpha}(0,T)$ the collection of process $(Y,Z,A)$ such that $Y\in S_G^{\alpha}(0,T)$, $Z\in H_G^{\alpha}(0,T)$, $A$ is a continuous nondecreasing process with $A_0=0$ and $A\in S_G^\alpha(0,T)$.  Note that we use a ``martingale condition" (c) instead of the Skorohod condition. Under some appropriate assumptions, we can prove that the solution of the above reflected $G$-BSDE is unique. In proving the existence of this problem, we should use the approximation method via penalization. This is a constructive method in the sense that the solution of the reflected $G$-BSDE is proved to be the limit of a sequence of penalized $G$-BSDEs. Different from the classical case, the dominated convergence theorem does not hold under $G$-framework. Besides, any bounded sequence in $M_G^p(0,T)$ is no longer weakly compact. The main difficulty in carrying out this construction is to prove the convergence property in some appropriate sense. It is worth pointing out that the main idea is to apply the uniformly continuous property of the elements in $S_G^p(0,T)$. 

Actually, the above equations hold $P$-$a.s.$ for every probability measure $P$ belongs to a non-dominated class of mutually singular measures. Therefore, the $G$-expectation theory shares many similarities with second order BSDEs (2BSDEs for short) developed by Cheridito, Soner, Touzi and Victoir \cite{CSTV}.  Matoussi, Possamai and Zhou \cite{MPZ} showed the existence and uniqueness of second order reflected BSDE whose solution is $(Y,Z,K^P)_{P\in\mathcal{P}_H^\kappa}$ satisfying
\[Y_t=\xi+\int_t^T \hat{F}_s(Y_s,Z_s)ds-\int_t^T Z_sdB_s+(K_T^P-K_t^P),\ P\textrm{-}a.s.,\]
with
\[Y_t\geq S_t, \ K_t^P-k_t^P=\underset{P'\in\mathcal{P}_H(t+,P)}{{\essinf}^P} E_t^{P'}[K_T^P-k_T^P], \ P\textrm{-}a.s., \ 0\leq t\leq T, \ \forall P\in\mathcal{P}_H^\kappa,\]
where $(y^P,z^p,k^P)$ denotes the unique solution to the standard RBSDE with data $(\xi,\hat{F},S)$ under $P$. The main contribution of our paper is that the triple $(Y,Z,A)$ is universally defined within the $G$-framework such that the processes have strong regularity property. Due to this property, the solution is time-consistent and the process $A$ can be aggregated into a universal process.

Similar with \cite{KKPPQ}, when the reflected $G$-BSDE is formulated under a Markovian framework, the solution of this problem provides a probabilistic representation for the solution of an obstacle problem for nonlinear parabolic PDE. There has been tremendous interest in developing the obstacle problem for partial differential equations since it has wide applications to mathematical finance (see \cite{EPQ}) and mathematical physics (see \cite{R87}). The method in this paper is called the Feynman-Kac formula which gives a link between probability theory and PDEs using the language of viscosity solutions. The other approach is related to the variational inequalities and in this case the solutions belong to a Sobolev space, see in particular \cite{KS} and \cite{F}. We should point out that both of the solutions studied by these two methods are in weak sense.

The rest of paper is organized as follows. In Section 2, we present some notations and results as preliminaries for the later proofs. The problem is formulated in detail in  Section 3 and we state some estimates of the solutions from which we derive some integrability properties of the solutions. In Section 4, we establish the approximation method via penalization. We state some convergence properties of the solution to the penalized $G$-BSDE. Our main results are showed and proved in Section 5. Furthermore, we prove a comparison theorem similar to that in \cite{HJPS2}, specifically for nonreflected $G$-BSDEs. In Section 6, we give the relation between reflected $G$-BSDEs and the corresponding obstacle problems for fully nonlinear parabolic PDEs. Finally, we use the results of the previous section to study the pricing problem for American contingent claims under model uncertainty in Section 7. In Appendix we introduce the optional stopping theorem under $G$-framework using for pricing for American contingent claims.

\section{Preliminaries}
We recall some basic notions and results of $G$-expectation, which are needed in the sequel. More relevant details can be found in \cite{HJPS1}, \cite{HJPS2},
\cite{P07a}, \cite{P08a}, \cite{P10}.

\subsection{$G$-expectation}

\begin{definition}
\label{def2.1} Let $\Omega$ be a given set and let $\mathcal{H}$ be a vector
lattice of real valued functions defined on $\Omega$, namely $c\in \mathcal{H%
}$ for each constant $c$ and $|X|\in \mathcal{H}$ if $X\in \mathcal{H}$. $%
\mathcal{H}$ is considered as the space of random variables. A sublinear
expectation $\hat{\mathbb{E}}$ on $\mathcal{H}$ is a functional $
\hat {\mathbb{E}}:\mathcal{H}\rightarrow \mathbb{R}$ satisfying the following
properties: for all $X,Y\in \mathcal{H}$, we have

\begin{description}
\item[(i)] Monotonicity: If $X\geq Y$, then $\hat{\mathbb{E}}[X]\geq
\hat{\mathbb{E}}[Y]$;

\item[(ii)] Constant preserving: $\hat{\mathbb{E}}[c]=c$;

\item[(iii)] Sub-additivity: $\hat{\mathbb{E}}[X+Y]\leq \hat{\mathbb{E}}[X]+%
\hat{\mathbb{E}}[Y]$;

\item[(iv)] Positive homogeneity: $\hat{\mathbb{E}}[\lambda X]=\lambda
\hat{\mathbb{E}}[X]$ for each $\lambda \geq0$.
\end{description}

The triple $(\Omega,\mathcal{H},\hat{\mathbb{E}})$ is called a
sublinear expectation space.  $X \in\mathcal{ H}$ is called a random
variable in $(\Omega,\mathcal{H},\hat{\mathbb{E}})$. We often call
$Y = (Y_1, \ldots, Y_d), Y_i \in\mathcal{ H}$ a $d$-dimensional
random vector in $(\Omega,\mathcal{H},\hat{\mathbb{E}})$.
\end{definition}

Let $\Omega_T=C_{0}([0,T];\mathbb{R}^{d})$, the space of
$\mathbb{R}^{d}$-valued continuous functions on $[0,T]$ with $\omega_{0}=0$, be endowed
with the supremum norm, and $B=(B^i)_{i=1}^d$ be the canonical
process. For each $T>0$, denote
\[
L_{ip} (\Omega_T):=\{ \varphi(B_{t_{1}},...,B_{t_{n}}):n\geq1,t_{1}%
,...,t_{n}\in\lbrack0,T],\varphi\in C_{Lip}(\mathbb{R}^{d\times n})\}.
\]

Denote by $\mathbb{S}_d$ the collection of all $d\times d$ symmetric matrices. For each given monotonic and sublinear function
$G:\mathbb{S}_{d}\rightarrow\mathbb{R}$, we can construct a $G$-expectation space  $(\Omega_T,  L_{ip}(\Omega_T),\hat{\mathbb{E}},\hat{\mathbb{E}}_t)$. The canonical process $B$ is the d-dimensional $G$-Brownian motion under this space. In this paper, we suppose that $G$ is non-degenerate, i.e., there exists some $\underline{\sigma}^2>0$ such that $G(A)-G(B)\geq \frac{1}{2}\underline{\sigma}^2 tr[A-B]$ for any $A\geq B$.

Let $B$ be the d-dimensional $G$-Brownian motion. For each fixed $a\in \mathbb{R}^d$,  $\{B^a_t\}:=\{\langle a, B_t\rangle\}$ is a 1-dimensional $G_a$-Brownian motion, where $G_a:\mathbb{R}\rightarrow \mathbb{R}$ satisfies
\[G_a(p)=G(aa^T)p^++G(-aa^T)p^-.\]
Let $\pi_{t}^{N}=\{t_{0}^{N},\cdots,t_{N}^{N}%
\}$, $N=1,2,\cdots$, be a sequence of partitions of $[0,t]$ such that $\mu
(\pi_{t}^{N})=\max\{|t_{i+1}^{N}-t_{i}^{N}|:i=0,\cdots,N-1\} \rightarrow0$,
the quadratic variation process of $B^a$ is defined by%
\[
\langle B^a\rangle_{t}=\lim_{\mu(\pi_{t}^{N})\rightarrow0}%
\sum_{j=0}^{N-1}(B^a_{t_{j+1}^{N}}-B^a_{t_{j}^{N}}%
)^{2}.
\]
For $a,\bar{a}\in\mathbb{R}^d$, we can define the mutual variation process of $B^a$ and $B^{\bar{a}}$ by
\[\langle B^a, B^{\bar{a}}\rangle_t:=\frac{1}{4}[\langle B^{a+\bar{a}}\rangle-\langle B^{a-\bar{a}}\rangle].\]

Denote by $L_{G}^{p}(\Omega_T)$   the completion of
$L_{ip} (\Omega_T)$ under the norm $\Vert\xi\Vert_{L_{G}^{p}}:=(\hat{\mathbb{E}}[|\xi|^{p}])^{1/p}$ for $p\geq1$. For all $\ t\in\lbrack
0,T]$, $\mathbb{\hat{E}}_{t}[\cdot]$ is a continuous mapping on $L_{ip}(\Omega_T)$
w.r.t. the norm $\Vert\cdot\Vert_{L_G^p}$. Therefore it can be
extended continuously to the completion $L_{G}^{p}(\Omega_{T})$. Denis et al. \cite{DHP11} proved the following representation theorem of $G$-expectation on $L_G^1(\Omega_T)$.
\begin{theorem}[\cite{DHP11,HP09}]
\label{the1.1}  There exists a weakly compact set
$\mathcal{P}\subset\mathcal{M}_{1}(\Omega_T)$, the set of all probability
measures on $(\Omega_T,\mathcal{B}(\Omega_T))$, such that
\[
\hat{\mathbb{E}}[\xi]=\sup_{P\in\mathcal{P}}E_{P}[\xi]\ \ \text{for
\ all}\ \xi\in  {L}_{G}^{1}{(\Omega_T)}.
\]
$\mathcal{P}$ is called a set that represents $\hat{\mathbb{E}}$.
\end{theorem}

Let $\mathcal{P}$ be a weakly compact set that represents $\hat{\mathbb{E}}$.
For this $\mathcal{P}$, we define capacity%
\[
c(A):=\sup_{P\in\mathcal{P}}P(A),\ A\in\mathcal{B}(\Omega_T).
\]

\begin{definition}A set $A\subset\mathcal{B}(\Omega_T)$ is polar if $c(A)=0$.  A
property holds $``quasi$-$surely"$ (q.s.) if it holds outside a
polar set.
\end{definition}

 In the following, we do not distinguish two random variables $X$ and $Y$ if $X=Y$ q.s..

For $\xi\in L_{ip}(\Omega_T)$, let $\mathcal{E}(\xi)=\hat{\mathbb{E}}[\sup_{t\in[0,T]}\hat{\mathbb{E}}_t[\xi]]$, where $\hat{\mathbb{E}}$ is the $G$-expectation. For convenience, we call $\mathcal{E}$ $G$-evaluation. For $p\geq 1$ and $\xi\in L_{ip}(\Omega_T)$, define $\|\xi\|_{p,\mathcal{E}}=[\mathcal{E}(|\xi|^p)]^{1/p}$ and denote by $L_{\mathcal{E}}^p(\Omega_T)$ the completion of $L_{ip}(\Omega_T)$ under $\|\cdot\|_{p,\mathcal{E}}$. The following estimate between the two norms $\|\cdot\|_{L_G^p}$ and $\|\cdot\|_{p,\mathcal{E}}$ will be frequently used in this paper.


\begin{theorem}[\cite{S11}]\label{the1.2}
For any $\alpha\geq 1$ and $\delta>0$, $L_G^{\alpha+\delta}(\Omega_T)\subset L_{\mathcal{E}}^{\alpha}(\Omega_T)$. More precisely, for any $1<\gamma<\beta:=(\alpha+\delta)/\alpha$, $\gamma\leq 2$, we have
\begin{displaymath}
\|\xi\|_{\alpha,\mathcal{E}}^{\alpha}\leq \gamma^*\{\|\xi\|_{L_G^{\alpha+\delta}}^{\alpha}+14^{1/\gamma}
C_{\beta/\gamma}\|\xi\|_{L_G^{\alpha+\delta}}^{(\alpha+\delta)/\gamma}\},\quad \forall \xi\in L_{ip}(\Omega_T).
\end{displaymath}
where $C_{\beta/\gamma}=\sum_{i=1}^\infty i^{-\beta/\gamma}$,$\gamma^*=\gamma/(\gamma-1)$.
\end{theorem}

\subsection{$G$-It\^{o} calculus}
\begin{definition}
\label{def2.6} Let $M_{G}^{0}(0,T)$ be the collection of processes in the
following form: for a given partition $\{t_{0},\cdot\cdot\cdot,t_{N}\}=\pi
_{T}$ of $[0,T]$,
\[
\eta_{t}(\omega)=\sum_{j=0}^{N-1}\xi_{j}(\omega)\mathbf{1}_{[t_{j},t_{j+1})}(t),
\]
where $\xi_{i}\in L_{ip}(\Omega_{t_{i}})$, $i=0,1,2,\cdot\cdot\cdot,N-1$. For each
$p\geq1$ and $\eta\in M_G^0(0,T)$ let $\|\eta\|_{H_G^p}:=\{\hat{\mathbb{E}}[(\int_0^T|\eta_s|^2ds)^{p/2}]\}^{1/p}$, $\Vert\eta\Vert_{M_{G}^{p}}:=(\mathbb{\hat{E}}[\int_{0}^{T}|\eta_{s}|^{p}ds])^{1/p}$ and denote by $H_G^p(0,T)$,  $M_{G}^{p}(0,T)$ the completion
of $M_{G}^{0}(0,T)$ under the norm $\|\cdot\|_{H_G^p}$, $\|\cdot\|_{M_G^p}$ respectively.
\end{definition}

For two processes $ \eta\in M_{G}^{2}(0,T)$ and $ \xi\in M_{G}^{1}(0,T)$,
the $G$-It\^{o} integrals  $(\int^{t}_0\eta_sdB^i_s)_{0\leq t\leq T}$ and $(\int^{t}_0\xi_sd\langle
B^i,B^j\rangle_s)_{0\leq t\leq T}$  are well defined, see  Li-Peng \cite{LP} and Peng \cite{P10}. Moreover, by Proposition 2.10 in \cite{LP} and the classical Burkholder-Davis-Gundy inequality, the following property holds.
\begin{proposition}\label{the1.3}
If $\eta\in H_G^{\alpha}(0,T)$ with $\alpha\geq 1$ and $p\in(0,\alpha]$, then we can get
$\sup_{u\in[t,T]}|\int_t^u\eta_s dB_s|^p\in L_G^1(\Omega_T)$ and
\begin{displaymath}
\underline{\sigma}^p c_p\hat{\mathbb{E}}_t[(\int_t^T |\eta_s|^2ds)^{p/2}]\leq
\hat{\mathbb{E}}_t[\sup_{u\in[t,T]}|\int_t^u\eta_s dB_s|^p]\leq
\bar{\sigma}^p C_p\hat{\mathbb{E}}_t[(\int_t^T |\eta_s|^2ds)^{p/2}].
\end{displaymath}
\end{proposition}

Let $S_G^0(0,T)=\{h(t,B_{t_1\wedge t}, \ldots,B_{t_n\wedge t}):t_1,\ldots,t_n\in[0,T],h\in C_{b,Lip}(\mathbb{R}^{n+1})\}$. For $p\geq 1$ and $\eta\in S_G^0(0,T)$, set $\|\eta\|_{S_G^p}=\{\hat{\mathbb{E}}[\sup_{t\in[0,T]}|\eta_t|^p]\}^{1/p}$. Denote by $S_G^p(0,T)$ the completion of $S_G^0(0,T)$ under the norm $\|\cdot\|_{S_G^p}$. We have the following continuity property for any $Y\in S_G^p(0,T)$ with $p>1$.

\begin{lemma}[\cite{LPS}]\label{the3.7}
For $Y\in S_G^p(0,T)$ with $p>1$, we have, by setting $Y_s:=Y_T$ for $s>T$,
\begin{displaymath}
F(Y):=\limsup_{\varepsilon\rightarrow0}(\hat{\mathbb{E}}[\sup_{t\in[0,T]}\sup_{s\in[t,t+\varepsilon]}|Y_t-Y_s|^p])^\frac{1}{p}=0.
\end{displaymath}
\end{lemma}

We now introduce some basic results of $G$-BSDEs. Consider the following type of $G$-BSDEs (here we use Einstein convention)
\begin{equation}\label{eq1.1}
Y_t=\xi+\int_t^T f(s,Y_s,Z_s)ds+\int_t^T g_{ij}(s,Y_s,Z_s)d\langle B^i,B^j\rangle_s-\int_t^T Z_s dB_s-(K_T-K_t),
\end{equation}
where
\begin{displaymath}
f(t,\omega,y,z),g_{ij}(t,\omega,y,z):[0,T]\times\Omega_T\times\mathbb{R}\times\mathbb{R}^d\rightarrow \mathbb{R},
\end{displaymath}
satisfying the following properties:
\begin{description}
\item[(H1')] There exists some $\beta>1$ such that for any $y,z$, $f(\cdot,\cdot,y,z),g_{ij}(\cdot,\cdot,y,z)\in M_G^{\beta}(0,T)$,
\item[(H2)] There exists some $L>0$ such that
\begin{displaymath}
|f(t,y,z)-f(t,y',z')|+\sum_{i,j=1}^d|g_{ij}(t,y,z)-g_{ij}(t,y',z')|\leq L(|y-y'|+|z-z'|).
\end{displaymath}
\end{description}

For simplicity, we denote by $\mathfrak{S}_G^{\alpha}(0,T)$ the collection of process $(Y,Z,K)$ such that $Y\in S_G^{\alpha}(0,T)$, $Z\in H_G^{\alpha}(0,T;\mathbb{R}^d)$, $K$ is a decreasing $G$-martingale with $K_0=0$ and $K_T\in L_G^{\alpha}(\Omega_T)$.


\begin{theorem}[\cite{HJPS1}]\label{the1.4}
Assume that $\xi\in L_G^{\beta}(\Omega_T)$ and $f,g_{ij}$ satisfy (H1') and (H2) for some $\beta>1$. Then for any $1<\alpha<\beta$, equation \eqref{eq1.1} has a unique solution $(Y,Z,K)\in \mathfrak{S}_G^{\alpha}(0,T)$.
\end{theorem}

We also have the comparison theorem for $G$-BSDE.
\begin{theorem}[\cite{HJPS2}]\label{the1.5}
Let $(Y_t^l,Z_t^l,K_t^l)_{t\leq T}$, $l=1,2$, be the solutions of the following $G$-BSDEs:
\begin{displaymath}
Y^l_t=\xi^l+\int_t^T f^l(s,Y^l_s,Z^l_s)ds+\int_t^T g^l_{ij}(s,Y^l_s,Z^l_s)d\langle B^i,B^j\rangle_s+V_T^l-V_t^l-\int_t^T Z^l_s dB_s-(K^l_T-K^l_t),
\end{displaymath}
where $\{V_t^l\}_{0\leq t\leq T}$ are RCLL processes such that $\hat{\mathbb{E}}[\sup_{t\in[0,T]}|V_t^l|^\beta]<\infty$, $f^l,\ g^l_{ij}$ satisfy (H1') and (H2), $\xi^l\in L_G^{\beta}(\Omega_T)$ with $\beta>1$. If $\xi^1\geq \xi^2$, $f^1\geq f^2$, $g^1_{ij}\geq g^2_{ij}$, for $i,j=1,\cdots,d$, $V_t^1-V_t^2$ is an increasing process, then $Y_t^1\geq Y_t^2$.
\end{theorem}


\section{Reflected $G$-BSDE with a lower obstacle and some a priori estimates}
For simplicity, we consider the $G$-expectation space $(\Omega,L_G^1(\Omega_T),\hat{\mathbb{E}})$ with $\Omega_T=C_0([0,T],\mathbb{R})$ and $\bar{\sigma}^2=\hat{\mathbb{E}}[B_1^2]\geq -\hat{\mathbb{E}}[-B_1^2]=\underline{\sigma}^2$. Our results and methods still hold for the case $d>1$. We are given the following data: the generator $f$ and $g$, the obstacle process $\{S_t\}_{t\in[0,T]}$ and the terminal value $\xi$, where $f$ and $g$ are maps
\begin{displaymath}
f(t,\omega,y,z),g(t,\omega,y,z):[0,T]\times\Omega_T\times\mathbb{R}^2\rightarrow\mathbb{R}.
\end{displaymath}

We will make the following assumptions: There exists some $\beta>2$ such that
\begin{description}
\item[(H1)] for any $y,z$, $f(\cdot,\cdot,y,z)$, $g(\cdot,\cdot,y,z)\in M_G^\beta(0,T)$;
\item[(H2)] $|f(t,\omega,y,z)-f(t,\omega,y',z')|+|g(t,\omega,y,z)-g(t,\omega,y',z')|\leq L(|y-y'|+|z-z'|)$ for some $L>0$;
\item[(H3)] $\xi\in L_G^\beta(\Omega_T)$ and $\xi\geq S_T$, $q.s.$;
\item[(H4)] There exists a constant $c$ such that $\{S_t\}_{t\in[0,T]}\in S_G^\beta(0,T)$ and $S_t\leq c$, for each $t\in[0,T]$;
\item[(H4')] $\{S_t\}_{t\in[0,T]}$ has the following form
\begin{displaymath}
S_t=S_0+\int_0^t b(s)ds+\int_0^t l(s)d\langle B\rangle_s+\int_0^t \sigma(s)dB_s,
\end{displaymath}
where $\{b(t)\}_{t\in[0,T]}$, $\{l(t)\}_{t\in[0,T]}$ belong to $M_G^\beta(0,T)$ and $\{\sigma(t)\}_{t\in[0,T]}$ belongs to $H_G^\beta(0,T)$.
\end{description}

Let us now introduce our reflected $G$-BSDE with a lower obstacle. A triplet of processes $(Y,Z,A)$ is called a solution of reflected $G$-BSDE with a lower obstacle if for some $1<\alpha\leq \beta$ the following properties hold:
\begin{description}
\item[(a)]$(Y,Z,A)\in\mathcal{S}_G^{\alpha}(0,T)$ and $Y_t\geq S_t$, $0\leq t\leq T$;
\item[(b)]$Y_t=\xi+\int_t^T f(s,Y_s,Z_s)ds+\int_t^T g(s,Y_s,Z_s)d\langle B\rangle_s
-\int_t^T Z_s dB_s+(A_T-A_t)$;
\item[(c)]$\{-\int_0^t (Y_s-S_s)dA_s\}_{t\in[0,T]}$ is a non-increasing $G$-martingale.
\end{description}
Here we denote by $\mathcal{S}_G^{\alpha}(0,T)$ the collection of process $(Y,Z,A)$ such that $Y\in S_G^{\alpha}(0,T)$, $Z\in H_G^{\alpha}(0,T;\mathbb{R})$, $A$ is a continuous nondecreasing process with $A_0=0$ and $A\in S_G^\alpha(0,T)$. For simplicity, we mainly consider the case with $g\equiv 0$ and $l\equiv 0$. Similar results still hold for the cases $g,l\neq 0$. Now we give a priori estimates for the solution of the reflected $G$-BSDE with a lower obstacle. In the following, $C$ will always designate a constant, which may vary from line to line.

\begin{proposition}\label{the1.6}
Let $f$ satisfies (H1) and (H2). Assume
\begin{displaymath}
Y_t=\xi+\int_t^T f(s,Y_s,Z_s)ds-\int_t^T Z_sdB_s+(A_T-A_t),
\end{displaymath}
where $Y\in S_G^{\alpha}(0,T)$, $Z\in H_G^{\alpha}(0,T;\mathbb{R})$, $A$ is a continuous nondecreasing process with $A_0=0$ and $A\in S_G^{\alpha}(0,T)$ for some $\alpha>1$. Then there exists a constant $C:=C(\alpha, T, L,\underline{\sigma})>0$ such that
\begin{align}
\hat{\mathbb{E}}_t[(\int_t^T |Z_s|^2ds)^{\frac{\alpha}{2}}]&\leq C\{\hat{\mathbb{E}}_t[\sup_{s\in[t,T]}|Y_s|^\alpha]
+(\hat{\mathbb{E}}_t[\sup_{s\in[t,T]}|Y_s|^\alpha])^{1/2}(\hat{\mathbb{E}}_t[(\int_t^T|f(s,0,0)|ds)^\alpha])^{1/2}\},\\
\hat{\mathbb{E}}_t[|A_T-A_t|^\alpha]&\leq C\{\hat{\mathbb{E}}_t[\sup_{s\in[t,T]}|Y_s|^\alpha]+\hat{\mathbb{E}}_t[(\int_t^T |f(s,0,0)|ds)^\alpha]\}.
\end{align}
\end{proposition}

\begin{proof}
The proof is similar to that of Proposition 3.5 in \cite{HJPS1}. So we omit it.
\end{proof}

\begin{proposition}\label{the1.7}
For $i=1,2$, let $\xi^i\in L_G^{\beta}(\Omega_T)$, $f^i$ satisfy (H1) and (H2) for some $\beta>2$. Assume
\begin{displaymath}
Y_t^i=\xi^i+\int_t^T f^i(s,Y_s,Z_s)ds-\int_t^T Z_s^idB_s+(A_T^i-A_t^i),
\end{displaymath}
where $Y^i\in S_G^\alpha(0,T)$, $Z^i\in H_G^\alpha(0,T)$, $A^i$ is a continuous nondecreasing process with $A_0^i=0$ and $A^i\in S_G^\alpha(0, T)$ for some $1<\alpha<\beta$. Set $\hat{Y}_t=Y^1_t-Y^2_t$, $\hat{Z}_t=Z^1_t-Z^2_t$, $\hat{A}_t=A^1_t-A^2_t$. Then there exists a constant $C:=C(\alpha,T,L,\underline{\sigma})$ such that
\begin{align*}
\hat{\mathbb{E}}[(\int_0^T|\hat{Z}|^2ds)^{\frac{\alpha}{2}}]\leq & C_{\alpha}\{(\hat{\mathbb{E}}[\sup_{t\in[0,T]}|\hat{Y}_t|^\alpha])^{1/2}
\sum_{i=1}^2[(\hat{\mathbb{E}}[\sup_{t\in[0,T]}|{Y}^i_t|^\alpha])^{1/2}\\&+(\hat{\mathbb{E}}[(\int_0^T |f^i(s,0,0)|ds)^\alpha])^{1/2}]+\hat{\mathbb{E}}[\sup_{t\in[0,T]}|\hat{Y}_t|^\alpha]\}.
\end{align*}
\end{proposition}
\begin{proof}
The proof is similar to that of Proposition 3.8 in \cite{HJPS1}. So we omit it.
\end{proof}
\begin{proposition}\label{the1.8}
For $i=1,2$, let $\xi^i\in L_G^{\beta}(\Omega_T)$ with $\xi^i\geq S_T^i$, where
\begin{displaymath}
S^i_t=S_0^i+\int_0^t b^i(s)ds+\int_0^t \sigma^i(s)dB_s.
\end{displaymath}
Here $\{b^i(s)\}\in M_G^\beta(0,T)$, $\{\sigma^i(s)\}\in H_G^\beta(0,T)$ for some $\beta>2$. Let $f^i$ satisfy (H1) and (H2). Assume that $(Y^i,Z^i,A^i)\in\mathcal{S}_G^\alpha(0,T)$ for some $1<\alpha<\beta$ are the solutions of the reflected $G$-BSDEs corresponding to $\xi^i$, $f^i$ and $S^i$. Set $\tilde{Y}_t=(Y^1_t-S_t^1)-(Y^2_t-S_t^2)$. Then there exists a constant $C:=C(\alpha,T,L,\underline{\sigma})$ such that
\begin{align*}
&|Y^{i}_t|^\alpha\leq C\hat{\mathbb{E}}_t[|\xi^{i}|^\alpha+\sup_{s\in[t,T]}|S^{i}_s|^\alpha +\int^T_t |\bar{\lambda}_s^{i,0}|^\alpha ds],\\
&|\tilde{Y}_t|^\alpha\leq C\hat{\mathbb{E}}_t[|\tilde{\xi}|^\alpha+\int_t^T(|\hat{\lambda}_s|^\alpha+|\hat{\rho}_s|^\alpha+|\hat{S}_{s}|^\alpha) ds],
\end{align*}
where $\tilde{\xi}=(\xi^{1}-S_T^1)-(\xi^{2}-S_T^2)$, $\hat{\lambda}_{s}=|f^{1}(s,Y_{s}^{2},Z_{s}%
^{2})-f^{2}(s,Y_{s}^{2},Z_{s}^{2})|$, $\hat{\rho}_s=|b^1(s)-b^2(s)|+|\sigma^1(s)-\sigma^2(s)|$, $\hat{S}_s=S^1_s-S^2_s$ and
$\bar{\lambda}_s^{i,0}=|f^{i}(s,0,0)|+|b^{i}(s)|+|\sigma^{i}(s)|.$
\end{proposition}

\begin{proof}
We will only show the second inequality, since the first one can be proved in a similar way.

For any $\varepsilon>0$, set $\hat{f}_t=f^1(t,Y_t^1,Z_t^1)-f^2(t,Y_t^2,Z_t^2)$, $\hat{f}^1_t=f^1(t,Y_t^1,Z_t^1)-f^1(t,Y_t^2,Z_t^2)$, $\hat{A}_t=A^1_t-A^2_t$, $\tilde{Z}_t=(Z^1_t-\sigma^1(t))-(Z^2_t-\sigma^2(t))$, $\varepsilon_\alpha=\varepsilon(1-\alpha/2)^+$ and $\bar{Y}_t=|\tilde{Y}_t|^2+\varepsilon_\alpha$. Applying It\^{o}'s formula to $\bar{Y}_t^{\frac{\alpha}{2}}e^{rt}$, where $r>0$ will be determined later, we get
\begin{equation}\label{0}
\begin{split}
&\quad \bar{Y}_t^{\alpha/2}e^{rt}+\int_t^T re^{rs}\bar{Y}_s^{\alpha/2}ds+\int_t^T \frac{\alpha}{2} e^{rs}
\bar{Y}_s^{\alpha/2-1}(\tilde{Z}_s)^2d\langle B\rangle_s\\
&=(\varepsilon_\alpha+|\tilde{\xi}|^2)^{\alpha/2}e^{rT}+
\alpha(1-\frac{\alpha}{2})\int_t^Te^{rs}\bar{Y}_s^{\alpha/2-2}(\tilde{Y}_s)^2(\tilde{Z}_s)^2d\langle B\rangle_s-\int_t^T\alpha e^{rs}\bar{Y}_s^{\alpha/2-1}\tilde{Y}_s\tilde{Z}_sdB_s\\
&\quad+\int_t^T{\alpha} e^{rs}\bar{Y}_s^{\alpha/2-1}\tilde{Y}_s(\hat{f}_s+b^1(s)-b^2(s))ds +\int_t^T\alpha e^{rs}\bar{Y}_s^{\alpha/2-1}\tilde{Y}_sd\hat{A}_s\\
&\leq(\varepsilon_\alpha+|\tilde{\xi}|^2)^{\alpha/2} e^{rT}+\int_t^T{\alpha}e^{rs}\bar{Y}_s^{\frac{\alpha-1}{2}}\{|\hat{f}_s^1+b^1(s)-b^2(s)|+\hat{\lambda}_s\}ds
\\
&\quad+\alpha(1-\frac{\alpha}{2})\int_t^Te^{rs}\bar{Y}_s^{\alpha/2-1}(\tilde{Z}_s)^2d\langle B\rangle_s-(M_T-M_t),
\end{split}
\end{equation}
where $M_t=\int_0^t \alpha e^{rs}\bar{Y}_s^{\alpha/2-1}(\tilde{Y}_s\tilde{Z}_sdB_s-(\tilde{Y}_s)^+dA_s^1-(\tilde{Y}_s)^-dA_s^2)$. We claim that $\{M_t\}$ is a $G$-martingale. Indeed, note that
\begin{displaymath}\tilde{Y}_t=Y^1_t-S^1_t+S_t^2-Y^2_t\leq Y^1_t-S_t^1.\end{displaymath}
Consequently, \begin{displaymath}(\tilde{Y}_t)^+\leq (Y^1_t-S_t^1)^+=Y^1_t-S_t^1.\end{displaymath}
Then we obtain \begin{displaymath}0\geq -\int_t^T (\tilde{Y}_s)^+ dA^1_s\geq -\int_t^T (Y^1_s-S_s^1)dA^1_s.\end{displaymath}
Thus we can conclude that
\begin{displaymath}
0\geq \hat{\mathbb{E}}_t[ -\int_t^T (\tilde{Y}_s)^+ dA^1_s]\geq \hat{\mathbb{E}}_t[-\int_t^T (Y^1_s-S_s^1)dA^1_s]=0.
\end{displaymath}
It follows that the process $\{K_t^1\}_{t\in[0,T]}=\{-\int_0^t (\tilde{Y}_s)^+ dA^1_s\}_{t\in[0,T]}$ is a non-increasing $G$-martingale. Set $\{K_t^2\}_{t\in[0,T]}=\{-\int_0^t (\tilde{Y}_s)^- dA^2_s\}_{t\in[0,T]}$. Both $\{K_t^1\}$ and $\{K_t^2\}$ are non-increasing $G$-martingales, so is $\{\int_0^t \alpha e^{rs}\bar{Y}_s^{\alpha/2-1}(dK_s^1+dK_s^2)\}$, which yields that $\{M_t\}_{t\in[0,T]}$ is a $G$-martingale.
From the assumption of $f^1$, we derive that
\begin{equation}\label{1}
\begin{split}
&\quad\int_t^T{\alpha} e^{rs}\bar{Y}_s^{\frac{\alpha-1}{2}}|\hat{f}_s^1+b^1(s)-b^2(s)|ds\\
&\leq \int_t^T{\alpha} e^{rs}\bar{Y}_s^{\frac{\alpha-1}{2}}\{L(|\tilde{Y}_s|+|\tilde{Z}_s|)+(L\vee 1)(|\hat{S}_s|+|\hat{\rho}_s|)\}ds\\
&\leq (\alpha L+\frac{\alpha L^2}{\underline{\sigma}^2(\alpha-1)})\int_t^T e^{rs}\bar{Y}_s^{\alpha/2}ds
+\frac{\alpha(\alpha-1)}{4}\int_t^Te^{rs}\bar{Y}_s^{\alpha/2-1}(\tilde{Z}_s)^2d\langle B\rangle_s\\
&\quad+(L\vee 1)\int_t^T{\alpha} e^{rs}\bar{Y}_s^{\frac{\alpha-1}{2}}\{|\hat{S}_s|+|\hat{\rho}_s|\}ds.
\end{split}
\end{equation}
By Young's inequality, we have
\begin{equation}\label{2}
\begin{split}
&\int_t^T{\alpha} e^{rs}\bar{Y}_s^{\frac{\alpha-1}{2}}\{|\hat{\lambda}_s|+|\hat{S}_s|+|\hat{\rho}_s|\}ds\\
 \leq &3(\alpha-1)\int_t^T  e^{rs}\bar{Y}_s^{\alpha/2}ds+\int_t^T e^{rs}\{|\hat{\lambda}_s|^\alpha+|\hat{\rho}_s|^\alpha+|\hat{S}_{s}|^\alpha\}ds.
\end{split}
\end{equation}
By \eqref{0}-\eqref{2} and setting $r=3(L\vee 1)(\alpha-1)+\alpha L+\frac{\alpha L^2}{\underline{\sigma}^2(\alpha-1)}+1$, we get
\begin{displaymath}
\bar{Y}_t^{\alpha/2}e^{rt}+(M_T-M_t)\leq C\{(\varepsilon_\alpha+|\tilde{\xi}|^2)^{\alpha/2} e^{rT}+\int_t^Te^{rs}(|\hat{\lambda}_s|^\alpha+|\hat{\rho}_s|^\alpha+|\hat{S}_{s}|^\alpha) ds\}.
\end{displaymath}
Taking conditional expectation on both sides and then by letting $\varepsilon\downarrow 0$, we have
\begin{align*}
|\tilde{Y}_t|^\alpha\leq C\hat{\mathbb{E}}_t[|\tilde{\xi}|^\alpha+\int_t^T(|\hat{\lambda}_s|^\alpha+|\hat{\rho}_s|^\alpha+|\hat{S}_{s}|^\alpha) ds].
\end{align*}
The proof is complete.
\end{proof}

\begin{proposition}\label{the1.9}
Let $(\xi,f,S)$ satisfy (H1)-(H4). 
Assume that $(Y,Z,A)\in\mathcal{S}_G^\alpha(0,T)$, for some $2\leq\alpha<\beta$, is a solution of the reflected $G$-BSDE with data $(\xi,f,S)$. Then there exists a constant $C:=C(\alpha,T, L,\underline{\sigma},c)>0$ such that
\begin{displaymath}
|Y_t|^\alpha\leq C\hat{\mathbb{E}}_t[1+|\xi|^\alpha+\int_t^T|f(s,0,0)|^\alpha ds].
\end{displaymath}
\end{proposition}

\begin{proof}
For any $r>0$, set $\tilde{Y}_t=|Y_t-c|^2$.  Applying It\^{o}'s formula to $\tilde{Y}_t^{\alpha/2} e^{rt}$, noting that $S_t\leq c$ and $A$ is a nondecreasing process, we have
\begin{align*}
&\tilde{Y}_t^{\alpha/2}e^{rt}+\int_t^T r e^{rs}\tilde{Y}_s^{\alpha/2}ds+\frac{\alpha}{2}\int_t^T e^{rs}\tilde{Y}_s^{\alpha/2-1} Z_s^2d\langle B\rangle_s\\
=&|\xi-c|^{\alpha}e^{rT}+\int_t^T \alpha e^{rs}\tilde{Y}_s^{\alpha/2-1}(Y_s-c) f(s,Y_s,Z_s)ds+\alpha(1-\frac{\alpha}{2})\int_t^T e^{rs}\tilde{Y}_s^{\alpha/2-2}(Y_s-c)^2Z_s^2\langle B\rangle_s\\
&-\int_t^T \alpha e^{rs}\tilde{Y}_s^{\alpha/2-1}(Y_s-c)Z_sdB_s+\int_t^T \alpha e^{rs}\tilde{Y}_s^{\alpha/2-1}(Y_s-c)dA_s\\
\leq &|\xi-c|^\alpha e^{rT}+\int_t^T \alpha e^{rs}\tilde{Y}_s^{\frac{\alpha-1}{2}} |f(s,Y_s,Z_s)|ds+\alpha(1-\frac{\alpha}{2})\int_t^T e^{rs}\tilde{Y}_s^{\alpha/2-1}Z_s^2\langle B\rangle_s-(M_T-M_t),
\end{align*}
where $M_t=\int_t^T \alpha e^{rs}\tilde{Y}_s^{\alpha/2-1}(Y_s-c)Z_sdB_s-\int_t^T \alpha e^{rs}\tilde{Y}_s^{\alpha/2-1}(Y_s-S_s)dA_s$. By condition (c), $M$ is a $G$-martingale. From the assumption of $f$ and by the Young inequality, we get
\begin{equation}\label{3}\begin{split}
\int_t^T \alpha e^{rs}\tilde{Y}_s^{\frac{\alpha-1}{2}} |f(s,Y_s,Z_s)|ds
\leq &\int_t^T \alpha e^{rs}\tilde{Y}_s^{\frac{\alpha-1}{2}} [|f(s,c,0)|+L|\tilde{Y}_s|+L|Z_s|]ds\\
\leq &(\alpha L+\frac{\alpha L^2}{\underline{\sigma}^2(\alpha-1)})\int_t^T e^{rs}\tilde{Y}_s^{\alpha/2}ds+(\alpha-1)\int_t^T e^{rs}\tilde{Y}_s^{\alpha/2}ds\\
&+\frac{\alpha(\alpha-1)}{4}\int_t^T e^{rs}\tilde{Y}_s^{\alpha/2-1}Z_s^2\langle B\rangle_s+\int_t^T e^{rs}|f(s,c,0)|^\alpha ds.
\end{split}\end{equation}
Setting $r=\alpha+\alpha L+\frac{\alpha L^2}{\underline{\sigma}^2(\alpha-1)}$ and by the above analysis, we have
\begin{displaymath}
\tilde{Y}_t^{\alpha/2}e^{rt}+M_T-M_t\leq |\xi-c|^\alpha e^{rT}+\int_t^T e^{rs}|f(s,c,0)|^\alpha ds.
\end{displaymath}
Taking conditional expectations on both side yields that
\begin{displaymath}
|Y_t-c|^\alpha\leq C\hat{\mathbb{E}}_t[|\xi-c|^\alpha+\int_t^T|f(s,c,0)|^\alpha ds].
\end{displaymath}
Noting that for $p\geq1$, we have $|a+b|^p\leq 2^{p-1}(|a|^p+|b|^p)$. Then the proof is complete.
\end{proof}

\begin{proposition}\label{the1.10}
Let $(\xi^1,f^1,S^1)$ and $(\xi^2,f^2,S^2)$ be two sets of data, each one satisfying all the assumptions (H1)-(H4). Let $(Y^i,Z^i,A^i)\in\mathcal{S}_G^\alpha(0,T)$ be a solution of the reflected $G$-BSDE with data $(\xi^i,f^i,S^i)$, $i=1,2$ respectively with $2\leq \alpha<\beta$. Set $\hat{Y}_t=Y^1_t-Y^2_t$, $\hat{S}_t=S^1_t-S^2_t$, $\hat{\xi}=\xi^1-\xi^2$. Then there exists a constant $C:=C(\alpha,T, L,\underline{\sigma},c)>0$ such that
\begin{displaymath}
|\hat{Y}_t|^\alpha\leq C\{\hat{\mathbb{E}}_t[|\hat{\xi}|^\alpha+\int_t^T|\hat{\lambda}_s|^\alpha ds]
+(\hat{\mathbb{E}}_t[\sup_{s\in[t,T]}|\hat{S}_s|^\alpha])^{\frac{1}{\alpha}}\Psi_{t,T}^{\frac{\alpha-1}{\alpha}}\},
\end{displaymath}
where $\hat{\lambda}_s=|f^1(s,Y_s^2,Z_s^2)-f^2(s,Y_s^2,Z_s^2)|$ and
\begin{displaymath}
\Psi_{t,T}=\sum_{i=1}^2\hat{\mathbb{E}}_t[\sup_{s\in [t,T]}\hat{\mathbb{E}}_s[1+|\xi^i|^\alpha+\int_t^T |f^i(r,0,0)|^\alpha dr]].
\end{displaymath}
\end{proposition}
\begin{proof}
Set $\hat{Z}_t=Z_t^1-Z_t^2$, $\hat{f}_t=f^1(t,Y_t^1,Z_t^1)-f^2(t,Y_t^2,Z_t^2)$ and $\hat{f}^1_t=f^1(t,Y_t^1,Z_t^1)-f^1(t,Y_t^2,Z_t^2)$. For any $r>0$, by applying It\^{o}'s formula to $\bar{Y}_t^{\alpha/2}e^{rt}=(|\hat{Y}_t|^2)^{\alpha/2} e^{rt}$, we have
\begin{equation}\label{eq1.2}
\begin{split}
&\quad \bar{Y}_t^{\alpha/2}e^{rt}+\int_t^T re^{rs}\bar{Y}_s^{\alpha/2}ds+\int_t^T \frac{\alpha}{2} e^{rs}
\bar{Y}_s^{\alpha/2-1}(\hat{Z}_s)^2d\langle B\rangle_s\\
&=|\hat{\xi}|^\alpha e^{rT}+
\alpha(1-\frac{\alpha}{2})\int_t^Te^{rs}\bar{Y}_s^{\alpha/2-2}(\hat{Y}_s)^2(\hat{Z}_s)^2d\langle B\rangle_s-\int_t^T\alpha e^{rs}\bar{Y}_s^{\alpha/2-1}\hat{Y}_s\hat{Z}_sdB_s\\
&\quad+\int_t^T{\alpha} e^{rs}\bar{Y}_s^{\alpha/2-1}\hat{Y}_s\hat{f}_sds +\int_t^T\alpha e^{rs}\bar{Y}_s^{\alpha/2-1}\hat{Y}_sd\hat{A}_s\\
&\leq|\hat{\xi}|^\alpha e^{rT}+
\alpha(1-\frac{\alpha}{2})\int_t^Te^{rs}\bar{Y}_s^{\alpha/2-1}(\hat{Z}_s)^2d\langle B\rangle_s+\int_t^T\alpha e^{rs}\bar{Y}_s^{\alpha/2-1}\hat{S}_sd\hat{A}_s\\
&\quad+\int_t^T{\alpha}e^{rs}\bar{Y}_s^{\frac{\alpha-1}{2}}\{|\hat{f}_s^1|+|\hat{\lambda}_s|\}ds-(M_T-M_t),
\end{split}\end{equation}
where $M_t=\int_0^t\alpha e^{rs}\bar{Y}_s^{\alpha/2-1}\hat{Y}_s\hat{Z}_sdB_s -\int_0^t\alpha e^{rs}\bar{Y}_s^{\alpha/2-1}(\hat{Y}_s-\hat{S}_s)^-dA^2_s-\int_0^t\alpha e^{rs}\bar{Y}_s^{\alpha/2-1}(\hat{Y}_s-\hat{S}_s)^+dA^1_s$. By a similar analysis as the proof of Proposition \ref{the1.8}, we conclude that $\{M_t\}_{t\in[0,T]}$ is a $G$-martingale. By Young's inequality and the assumption of $f^1$, similar with inequalities \eqref{1} and \eqref{2}, we have
\begin{align*}
\int_t^T{\alpha}e^{rs}\bar{Y}_s^{\frac{\alpha-1}{2}}\{|\hat{f}_s^1|+|\hat{\lambda}_s|\}ds
\leq &\frac{\alpha(\alpha-1)}{4}\int_t^T e^{rs}\tilde{Y}_s^{\alpha/2-1}Z_s^2\langle B\rangle_s+\int_t^T e^{rs}|\hat{\lambda}_s|^\alpha ds\\
&+(\alpha-1+\alpha L+\frac{\alpha L^2}{\underline{\sigma}^2(\alpha-1)})\int_t^T e^{rs}\tilde{Y}_s^{\alpha/2}ds.
\end{align*}
Set $r=\alpha+\alpha L+\frac{\alpha L^2}{\underline{\sigma}^2(\alpha-1)}$. Taking conditional expectations on both sides of \eqref{eq1.2}, we obtain
\begin{displaymath}
|\hat{Y}_t|^\alpha\leq C\{\hat{\mathbb{E}}_t[|\hat{\xi}|^\alpha+\int_t^T|\hat{\lambda}_s|^\alpha ds]+\hat{\mathbb{E}}_t[\int_t^T \bar{Y}_s^{\alpha/2-1}|\hat{S}_s|d({A}^1_s+A_s^2)]\}.
\end{displaymath}
By applying H\"{o}lder's inequality, we get
\begin{align*}
\hat{\mathbb{E}}_t[\int_t^T \bar{Y}_s^{\alpha/2-1}|\hat{S}_s|d({A}^1_s+A_s^2)]&\leq \hat{\mathbb{E}}_t[\sup_{s\in[t,T]}\bar{Y}_s^{\alpha/2-1}|\hat{S}_s|(|A^1_T-A^1_t|+|A_T^2-A_t^2|)]\\
&\leq (\hat{\mathbb{E}}_t[\sup_{s\in[0,T]}|\hat{S}_s|^{\alpha}])^{\frac{1}{\alpha}}(\hat{\mathbb{E}}_t[\sup_{s\in[t,T]}\bar{Y}_s^{\alpha/2}])^{\frac{\alpha-2}{\alpha}}
(\sum_{i=1}^2\hat{\mathbb{E}}_t[|A_T^i-A_t^i|^\alpha])^{\frac{1}{\alpha}}.
\end{align*}
From Proposition \ref{the1.6} and Proposition \ref{the1.9}, we finally get the desired result.
\end{proof}

\begin{remark}
If we require that the solution of a reflected $G$-BSDE is a quadruple $\{(Y_t,Z_t,K_t,A_t),0\leq t\leq T\}$ satisfying conditions (1)-(3) in the introduction, the solution is not unique. We can see this fact from the following example.

Let $f\equiv -1$, $\xi=0$ and $S\equiv 0$. It is easy to check that $(0,0,0,t)$ and $(0,0,\frac{1}{\bar{\sigma}^2-\underline{\sigma}^2}(\underline{\sigma}^2t-\langle B\rangle_t),\frac{1}{\bar{\sigma}^2-\underline{\sigma}^2}(\bar{\sigma}^2t-\langle B\rangle_t))$ are solutions of reflected $G$-BSDE with data $(0,-1,0)$ satisfying all the conditions (1)-(3).
\end{remark}

\section{Penalized method and convergence properties}
In order to derive the existence of the solution to the reflected $G$-BSDE with a lower obstacle, we shall apply the approximation method via penalization. In this section, we first state some convergence properties of solutions to the penalized $G$-BSDEs, which will be needed in the sequel.

For $f$ and $\xi$ satisfy (H1)-(H3), $\{S_t\}_{t\in[0,T]}$ satisfies (H4) or (H4'), we now consider the following family of $G$-BSDEs parameterized by $n=1,2,\cdots$,
\begin{equation}\label{eq1.3}
Y_t^n=\xi+\int_t^T f(s,Y_s^n,Z_s^n)ds+n\int_t^T(Y_s^n-S_s)^-ds-\int_t^T Z_s^ndB_s-(K_T^n-K_t^n).
\end{equation}

Now let $L_t^n=n\int_0^t (Y_s^n-S_s)^-ds$, then $(L_t^n)_{t\in[0,T]}$ is a nondecreasing process. We can rewrite reflected $G$-BSDE \eqref{eq1.3} as
\begin{equation}\label{eq1.4}
Y_t^n=\xi+\int_t^T f(s,Y_s^n,Z_s^n)ds-\int_t^T Z_s^ndB_s-(K_T^n-K_t^n)+(L_T^n-L_t^n).
\end{equation}

We now establish a priori estimates on the sequence $(Y^n,Z^n,K^n,L^n)$ which are uniform in $n$.
\begin{lemma}\label{the1.11}
There exists a constant C independent of n, such that for $1<\alpha<\beta$,\begin{displaymath}
\hat{\mathbb{E}}[\sup_{t\in[0,T]}|Y_t^n|^\alpha]\leq C,\textrm{ } \hat{\mathbb{E}}[|K_T^n|^\alpha]\leq C,\textrm{ }\hat{\mathbb{E}}[|L_T^n|^\alpha]\leq C,
\textrm{ }\hat{\mathbb{E}}[(\int_0^T |Z_t^n|^2 dt)^{\frac{\alpha}{2}}]\leq C.
\end{displaymath}
\end{lemma}

\begin{proof}
For simplicity, first we consider the case $S\equiv 0$. The proof of the other cases will be given in the remark. $\forall r,\varepsilon>0$, set $\tilde{Y}_t=(Y_t^n)^2+\varepsilon_\alpha$,
where $\varepsilon_\alpha=\varepsilon(1-\alpha/2)^+$. Note that for any $a\in\mathbb{R}$, $a\times a^-\leq 0$. By applying It\^{o}' formula to $\tilde{Y}_t^{\alpha/2}e^{rt}$ yields that
\begin{displaymath}
\begin{split}
&\quad\tilde{Y}_t^{\alpha/2}e^{rt}+\int_t^T re^{rs}\tilde{Y}_s^{\alpha/2}ds+\int_t^T \frac{\alpha}{2} e^{rs}
\tilde{Y}_s^{\alpha/2-1}(Z_s^n)^2d\langle B\rangle_s\\
&=(|\xi|^2+\varepsilon_\alpha)^{\frac{\alpha}{2}}e^{rT}+
\alpha(1-\frac{\alpha}{2})\int_t^Te^{rs}\tilde{Y}_s^{\alpha/2-2}(Y_s^n)^2(Z_s^n)^2d\langle B\rangle_s
+\int_t^T\alpha e^{rs}\tilde{Y}_s^{\alpha/2-1}Y_s^ndL_s^n\\
&\quad+\int_t^T{\alpha} e^{rs}\tilde{Y}_s^{\alpha/2-1}Y_s^nf(s,Y_s^n,Z_s^n)ds-
\int_t^T\alpha e^{rs}\tilde{Y}_s^{\alpha/2-1}(Y_s^nZ_s^ndB_s+Y_s^ndK_s^n)\\
&\leq(|\xi|^2+\varepsilon_\alpha)^{\frac{\alpha}{2}}e^{rT}+
\alpha(1-\frac{\alpha}{2})\int_t^Te^{rs}\tilde{Y}_s^{\alpha/2-2}(Y_s^n)^2(Z_s^n)^2d\langle B\rangle_s\\
&\quad+\int_t^T{\alpha} e^{rs}\tilde{Y}_s^{\alpha/2-1/2}|f(s,Y_s^n,Z_s^n)|ds-(M_T-M_t),
\end{split}
\end{displaymath}
where $M_t=\int_t^T\alpha e^{rs}\tilde{Y}_s^{\alpha/2-1}(Y_s^nZ_s^ndB_s+(Y_s^n)^+dK_s^n)$ is a $G$-martingale.
 Similar with inequality \eqref{3}, we have
\begin{displaymath}
\begin{split}
\int_t^T{\alpha} e^{rs}\tilde{Y}_s^{\frac{\alpha-1}{2}}|f(s,Y_s^n,Z_s^n)|ds
\leq &\int_t^T e^{rs}|f(s,0,0)|^\alpha ds+\frac{\alpha(\alpha-1)}{4}\int_t^Te^{rs}\tilde{Y}_s^{\alpha/2-1}(Z_s^n)^2d\langle B\rangle_s\\
&+(\alpha-1+\alpha L+\frac{\alpha L^2}{\underline{\sigma}^2(\alpha-1)})\int_t^T e^{rs}\tilde{Y}_s^{\alpha/2}ds.
\end{split}
\end{displaymath}
Set $r=\alpha+\alpha L+\frac{\alpha L^2}{\underline{\sigma}^2(\alpha-1)}$. We derive that
\begin{displaymath}
\tilde{Y}_t^{\alpha/2}e^{rt}+M_T-M_t\leq (|\xi|^2+\varepsilon_\alpha)^{\frac{\alpha}{2}}e^{rT}+\int_t^T e^{rs}|f(s,0,0)|^\alpha ds.
\end{displaymath}
Taking conditional expectations on both sides and then by letting $\varepsilon\rightarrow 0$, we obtain
\begin{displaymath}
|Y_t^n|^\alpha\leq C\hat{\mathbb{E}}_t[|\xi|^\alpha+\int_t^T|f(s,0,0)|^\alpha ds].
\end{displaymath}
By Theorem \ref{the1.2}, for $1<\alpha<\beta$, there exists a constant $C$ independent of $n$ such that $\hat{\mathbb{E}}[\sup_{t\in[0,T]}|Y_t^n|^\alpha]\leq C$.
By Proposition \ref{the1.6}, we have
\begin{align*}
\hat{\mathbb{E}}[(\int_0^T |Z_s^n|^2ds)^{\frac{\alpha}{2}}]&\leq C_{\alpha}\{\hat{\mathbb{E}}[\sup_{t\in[0,T]}|Y_t^n|^\alpha]
+(\hat{\mathbb{E}}[\sup_{t\in[0,T]}|Y_t^n|^\alpha])^{\frac{1}{2}}(\hat{\mathbb{E}}[(\int_0^T|f(s,0,0)|ds)^\alpha])^{\frac{1}{2}}\},\\
\hat{\mathbb{E}}[|L_T^n-K_T^n|^\alpha]&\leq C_{\alpha}\{\hat{\mathbb{E}}[\sup_{t\in[0,T]}|Y_t^n|^\alpha]+\hat{\mathbb{E}}[(\int_0^T|f(s,0,0)|ds)^\alpha]\},
\end{align*}
where the constant $C_\alpha$ depends on $\alpha,T,\underline{\sigma}$ and $L$. Thus we conclude that there exists a constant $C$ independent of $n$, such that for $1<\alpha<\beta$,
\begin{displaymath}
\hat{\mathbb{E}}[(\int_0^T |Z_t^n|^2 dt)^{\frac{\alpha}{2}}]\leq C, \quad \hat{\mathbb{E}}[|L_T^n-K_T^n|^\alpha]\leq C.
\end{displaymath}
Since $L_T^n$ and $-K_T^n$ are nonnegative, it follows that
\begin{displaymath}
\hat{\mathbb{E}}[|K_T^n|^\alpha]\leq C, \quad \hat{\mathbb{E}}[|L_T^n|^\alpha]=n^\alpha\hat{\mathbb{E}}[(\int_0^T (Y_s^n)^-s)^\alpha]\leq C.
\end{displaymath}
\end{proof}

\begin{remark}\label{rem1.1}
If the obstacle process $\{S_t\}_{t\in[0,T]}$ satisfies (H4), set $\tilde{Y}_t^n={Y}_t^n-c$. It is simple to check that
\begin{displaymath}
\tilde{Y}_t^n=\xi-c+\int_t^T f(s,\tilde{Y}_s^n+c,Z_s^n)ds+\int_t^T n(\tilde{Y}_s^n-(S_s-c))^-ds-\int_t^T Z_s^ndB_s-(K_T^n-K_t^n).
\end{displaymath}
By a similar analysis as the proof of Lemma \ref{the1.11}, we derive that
\begin{displaymath}
|\tilde{Y}_t^n|^\alpha\leq C\hat{\mathbb{E}}_t[|\xi-c|^\alpha+\int_t^T |f(s,c,0)|^\alpha ds].
\end{displaymath}

If $S$ satisfies (H4'), for simplicity we suppose that $l\equiv 0$. Let $\widetilde{Y}^n_t={Y}_t^n-{S}_t$ and $\widetilde{Z}_t^n={Z}_t^n-\sigma(s)$, we can rewrite \eqref{eq1.3} as the following:
\begin{displaymath}
\widetilde{Y}_t^n=\xi-{S}_T+\int_t^T [f(s,\widetilde{Y}_s^n+{S}_s,\widetilde{Z}_s^n+\sigma(s))+b(s)]ds+n\int_t^T(\widetilde{Y}_s^n)^-ds-\int_t^T \widetilde{Z}_s^ndB_s-({K}_T^n-{K}_t^n).
\end{displaymath}
Using the same method, we get
\begin{displaymath}
|\widetilde{Y}_t^n|^\alpha\leq C\hat{\mathbb{E}}_t[|\xi-{S}_T|^\alpha+\int_t^T|f(s,{S}_s,\sigma(s))+b(s)|^\alpha ds].
\end{displaymath}

Thus we conclude that in the above two cases, for $1<\alpha<\beta$, there exists a constant $C$ independent of $n$ such that $\hat{\mathbb{E}}[\sup_{t\in[0,T]}|{Y}_t^n|^\alpha]\leq C$. By Proposition \ref{the1.6}, we have
\begin{displaymath}
\hat{\mathbb{E}}[|{K}_T^n|^\alpha]\leq C, \ \hat{\mathbb{E}}[|{L}_T^n|^\alpha]=n^\alpha\hat{\mathbb{E}}[(\int_0^T ({Y}_s^n-S_s)^-ds)^\alpha]\leq C, \textrm{ and }\hat{\mathbb{E}}[(\int_0^T |Z_t^n|^2 dt)^{\frac{\alpha}{2}}]\leq C.
\end{displaymath}

\end{remark}



Lemma \ref{the1.11} implies that $(Y^n-S)^-\rightarrow 0$ in $M_G^1(0,T)$.  The following lemma which corresponds to Lemma 6.1 in \cite{KKPPQ} shows that this convergence holds in $S_G^\alpha(0,T)$, for $1<\alpha<\beta$. It is of vital importance to prove the convergence property for $(Y^n)$.

\begin{lemma}\label{the1.12}
For some $1<\alpha<\beta$, we have
\begin{displaymath}
\lim_{n\rightarrow\infty}\hat{\mathbb{E}}[\sup_{t\in[0,T]}|({Y}^n_t-S_t)^-|^\alpha]=0.
\end{displaymath}
\end{lemma}

\begin{proof}
We now consider the following $G$-BSDEs parameterized by $n=1,2,\cdots$,
\begin{displaymath}
y_t^n=\xi+\int_t^T f(s,Y_s^n,Z_s^n)ds+\int_t^T n(S_s-y_s^n)ds-\int_t^T z_s^ndB_s-(k_T^n-k_t^n).
\end{displaymath}
By applying $G$-It\^{o}'s formula to $e^{-nt}y_t^n$, we can get
\begin{displaymath}
y_t^n=e^{nt}\hat{\mathbb{E}}_t[e^{-nT}\xi+\int_t^T ne^{-ns}S_s ds+\int_t^T e^{-ns}f(s,Y_s^n,Z_s^n)ds].
\end{displaymath}
By the comparison theorem \ref{the1.5}, we have for all $n\geq 1$, $Y_t^n\geq Y_t^1$ and
\begin{displaymath}
{Y}^n_t-S_t\geq y_t^n-S_t=\hat{\mathbb{E}}_t[\tilde{S}_t^n+\int_t^T e^{n(t-s)}f(s,Y_s^n,Z_s^n)ds],
\end{displaymath}
where $\tilde{S}_t^n=e^{n(t-T)}(\xi-S_t)+\int_t^T ne^{n(t-s)}(S_s-S_t )ds$. It follows that
\begin{displaymath}
({Y}^n_t-S_t)^-\leq (y_t^n-S_t)^- \leq \hat{\mathbb{E}}_t[|\tilde{S}_t^n|+|\int_t^T e^{n(t-s)}f(s,Y_s^n,Z_s^n)ds|].
\end{displaymath}
Applying H\"{o}lder's inequality yields that
\begin{align*}
|\int_t^T e^{n(t-s)}f(s,Y_s^n,Z_s^n)ds|&\leq \frac{1}{\sqrt{2n}}(\int_0^T f^2(s,Y_s^n,Z_s^n)ds)^{1/2}\\
&\leq \frac{C}{\sqrt{n}}(\sup_{s\in[0,T]}|Y_s^n|^2+\int_0^T f^2(s,0,0)+|Z_s^n|^2ds)^{1/2}.
\end{align*}
By Lemma \ref{the1.11}, for $1<\alpha<\beta$, we have
\begin{equation}\label{e2}
\hat{\mathbb{E}}[\sup_{t\in[0,T]}|\int_t^T e^{n(t-s)}f(s,Y_s^n,Z_s^n)ds|^\alpha]\rightarrow 0, \textrm{ as } n\rightarrow\infty.
\end{equation}
For $\varepsilon>0$, it is straightforward to show that
\begin{align*}
|\tilde{S}_t^n|&=|e^{n(t-T)}(\xi-S_t)+\int_{t+\varepsilon}^T ne^{n(t-s)}(S_s-S_t )ds+\int_t^{t+\varepsilon} ne^{n(t-s)}(S_s-S_t )ds|\\
&\leq e^{n(t-T)}|\xi-S_t|+e^{-n\varepsilon}\sup_{s\in[t+\varepsilon,T]}|S_t-S_s|+\sup_{s\in[t,t+\varepsilon]}|S_s-S_t|.
\end{align*}
For $T>\delta>0$, from the above inequality we obtain
\begin{align*}
\sup_{t\in[0,T-\delta]}|\tilde{S}_t^n|&\leq e^{-n\delta}\sup_{t\in[0,T-\delta]}|\xi-S_t|+e^{-n\varepsilon}\sup_{t\in[0,T-\delta]}\sup_{s\in[t+\varepsilon,T]}|S_t-S_s|+\sup_{t\in[0,T-\delta]}\sup_{s\in[t,t+\varepsilon]}|S_s-S_t|\\
&\leq e^{-n\delta}(\sup_{t\in[0,T]}|S_t|+|\xi|)+2e^{-n\varepsilon}\sup_{t\in[0,T]}|S_t|+\sup_{t\in[0,T]}\sup_{s\in[t,t+\varepsilon]}|S_s-S_t|.
\end{align*}
It is easy to check that for each fixed $\varepsilon,\delta>0$,
\begin{equation}\begin{split}\label{e3}
\hat{\mathbb{E}}[\sup_{t\in[0,T-\delta]}|\tilde{S}_t^n|^\beta]&\leq C\{(e^{-n\beta\varepsilon}+e^{-n\beta\delta})\hat{\mathbb{E}}[\sup_{t\in[0,T]}|S_t|^\beta+|\xi|^\beta]+\hat{\mathbb{E}}[\sup_{t\in[0,T]}\sup_{s\in[t,t+\varepsilon]}|S_s-S_t|^\beta]\}\\
&\rightarrow C\hat{\mathbb{E}}[\sup_{t\in[0,T]}\sup_{s\in[t,t+\varepsilon]}|S_s-S_t|^\beta], \quad \textrm{as } n\rightarrow \infty.
\end{split}\end{equation}
For $1<\alpha<\beta$ and $0<\delta<T$, we have
\begin{equation}\begin{split}\label{e4}
& \hat{\mathbb{E}}[\sup_{t\in[0,T]}|({Y}_t^n-S_t)^-|^\alpha]
\leq \hat{\mathbb{E}}[\sup_{t\in[0,T-\delta]}|(Y_t^n-S_t)^-|^\alpha]+\hat{\mathbb{E}}[\sup_{t\in[T-\delta,T]}|(Y_t^n-S_t)^-|^\alpha]\\
\leq &\hat{\mathbb{E}}[\sup_{t\in[0,T-\delta]}\{\hat{\mathbb{E}}_t[|\tilde{S}_t^n|+|\int_t^T e^{n(t-s)}f(s,Y_s^n,Z_s^n)ds|]\}^\alpha]+\hat{\mathbb{E}}[\sup_{t\in[T-\delta,T]}|(Y_t^1-S_t)^-|^\alpha]\\
\leq & C\{\hat{\mathbb{E}}[\sup_{t\in[0,T-\delta]}\hat{\mathbb{E}}_t[\sup_{u\in[0,T-\delta]}|\tilde{S}_u^n|^\alpha]]+\hat{\mathbb{E}}[\sup_{t\in[0,T-\delta]}\hat{\mathbb{E}}_t[\sup_{u\in[0,T]}|\int_u^T e^{n(t-s)}f(s,Y_s^n,Z_s^n)ds|^\alpha]]\}\\
&+\hat{\mathbb{E}}[\sup_{t\in[T-\delta,T]}|(Y_t^1-S_t)^-|^\alpha]=:I+\hat{\mathbb{E}}[\sup_{t\in[T-\delta,T]}|(Y_t^1-S_t)^-|^\alpha].
\end{split}\end{equation}
By Lemma \ref{the3.7}, noting that $Y^1-S\in S_G^\alpha(0,T)$ and $(Y^1_T-S_T)^-=0$, we obtain
\begin{displaymath}
\lim_{\delta\rightarrow 0}\hat{\mathbb{E}}[\sup_{t\in[T-\delta,T]}|(Y_t^1-S_t)^-|^\alpha]=0.
\end{displaymath}
By Theorem \ref{the1.2} and combining \eqref{e2}, \eqref{e3}, we derive that
\begin{displaymath}
I\leq C\{\hat{\mathbb{E}}[\sup_{t\in[0,T]}\sup_{s\in[t,t+\varepsilon]}|S_s-S_t|^\beta]+(\hat{\mathbb{E}}[\sup_{t\in[0,T]}\sup_{s\in[t,t+\varepsilon]}|S_s-S_t|^\beta])^{\alpha/\beta}\}, \textrm{ as } n\rightarrow \infty.
\end{displaymath}
Now first let $n\rightarrow\infty$ and then let $\varepsilon,\delta\rightarrow0$ in \eqref{e4}. By Lemma \ref{the3.7} again, the above analysis proves that for $1<\alpha<\beta$,
\begin{displaymath}\lim_{n\rightarrow\infty}\hat{\mathbb{E}}[\sup_{t\in[0,T]}|({Y}_t^n-S_t)^-|^\alpha]=0.
\end{displaymath}\end{proof}

Now we show the convergence property of sequence $(Y^n)_{n=1}^{\infty}$.
\begin{lemma}\label{the1.13}
For some $\beta>\alpha\geq2$, we have
\begin{displaymath}
\lim_{n,m\rightarrow\infty}\hat{\mathbb{E}}[\sup_{t\in[0,T]}|Y_t^n-Y_t^m|^\alpha]=0.
\end{displaymath}
\end{lemma}

\begin{proof}
Without loss of generality, we may assume $S\equiv 0$ in \eqref{eq1.3}. For any $r>0$, set $\bar{Y}_t=|Y^n_t-Y^m_t|^2$, $\hat{f}_t=f(t,Y_t^n,Z_t^n)-f(t,Y_t^m,Z_t^m)$. By applying It\^{o}'s formula to $\bar{Y}_t^{\alpha/2}e^{rt}$, we  get
\begin{displaymath}
\begin{split}
&\quad\bar{Y}_t^{\alpha/2}e^{rt}+\int_t^T re^{rs}\bar{Y}_s^{\alpha/2}ds+\int_t^T \frac{\alpha}{2} e^{rs}
\bar{Y}_s^{\alpha/2-1}(\hat{Z}_s)^2d\langle B\rangle_s\\
&=
\alpha(1-\frac{\alpha}{2})\int_t^Te^{rs}\bar{Y}_s^{\alpha/2-2}(\hat{Y}_s)^2(\hat{Z}_s)^2d\langle B\rangle_s
+\int_t^T\alpha e^{rs}\bar{Y}_s^{\alpha/2-1}\hat{Y}_sd\hat{L}_s\\
&\quad+\int_t^T{\alpha} e^{rs}\bar{Y}_s^{\alpha/2-1}\hat{Y}_s\hat{f}_sds-\int_t^T\alpha e^{rs}\bar{Y}_s^{\alpha/2-1}(\hat{Y}_s\hat{Z}_sdB_s+\hat{Y}_sd\hat{K}_s)\\
&\leq
\alpha(1-\frac{\alpha}{2})\int_t^Te^{rs}\bar{Y}_s^{\alpha/2-2}(\hat{Y}_s)^2(\hat{Z}_s)^2d\langle B\rangle_s+\int_t^T{\alpha} e^{rs}\bar{Y}_s^{\frac{\alpha-1}{2}}|\hat{f}_s|ds\\
&\quad-\int_t^T\alpha e^{rs}\bar{Y}_s^{\alpha/2-1}Y_s^ndL^m_s-\int_t^T\alpha e^{rs}\bar{Y}_s^{\alpha/2-1}Y_s^mdL^n_s-(M_T-M_t),
\end{split}
\end{displaymath}
where $M_t=\int_0^t \alpha e^{rs}\bar{Y}_s^{\alpha/2-1}(\hat{Y}_s\hat{Z}_sdB_s+(\hat{Y}_s)^+dK_s^m+(\hat{Y}_s)^-dK_s^n)$ is a $G$-martingale. Similar with \eqref{1}, we have
\begin{displaymath}
\int_t^T{\alpha} e^{rs}\bar{Y}_s^{\frac{\alpha-1}{2}}|\hat{f}_s|ds
\leq (\alpha L+\frac{\alpha L^2}{\underline{\sigma}^2(\alpha-1)})\int_t^T e^{rs}\bar{Y}_s^{\alpha/2}ds
+\frac{\alpha(\alpha-1)}{4}\int_t^Te^{rs}\bar{Y}_s^{\alpha/2-1}(\hat{Z}_s)^2d\langle B\rangle_s.
\end{displaymath}
Let $r=1+\alpha L+\frac{\alpha L^2}{\underline{\sigma}^2(\alpha-1)}$. By the above analysis, we have
\[\bar{Y}_t^{\alpha/2}e^{rt}+(M_T-M_t)\leq-\int_t^T\alpha e^{rs}\bar{Y}_s^{\alpha/2-1}Y_s^ndL^m_s-\int_t^T\alpha e^{rs}\bar{Y}_s^{\alpha/2-1}Y_s^mdL^n_s.\]
Then taking conditional expectation on both sides of the above inequality, we conclude that
\begin{equation}\label{eq1.5}
\bar{Y}_t^{\alpha/2}e^{rt}\leq
\hat{\mathbb{E}}_t[-\int_t^T\alpha e^{rs}\bar{Y}_s^{\alpha/2-1}Y_s^ndL^m_s-\int_t^T\alpha e^{rs}\bar{Y}_s^{\alpha/2-1}Y_s^mdL^n_s].
\end{equation}
Observe that
\begin{displaymath}
\begin{split}
&\quad\hat{\mathbb{E}}_t[-\int_t^T\alpha e^{rs}\bar{Y}_s^{\alpha/2-1}Y_s^mdL^n_s]
\leq \alpha e^{rT}\hat{\mathbb{E}}_t[\int_t^T\bar{Y}_s^{\alpha/2-1}n(Y_s^n)^-(Y_s^m)^-ds]\\
&\leq C\hat{\mathbb{E}}_t[\int_0^Tn|(Y_s^n)^-|^{\alpha-1}(Y_s^m)^-ds]+
C\hat{\mathbb{E}}_t[\int_0^Tn|(Y_s^m)^-|^{\alpha-1}(Y_s^n)^-ds].
\end{split}
\end{displaymath}
From \eqref{eq1.5} and taking expectations on both sides, we deduce that
\begin{equation}\begin{split}\label{eq1.6}
\hat{\mathbb{E}}[\sup_{t\in[0,T]}|Y_t^n-Y_t^m|^\alpha]\leq
&C\hat{\mathbb{E}}[\sup_{t\in[0,T]}\{\hat{\mathbb{E}}_t[\int_0^T(n+m)|(Y_s^n)^-|^{\alpha-1}(Y_s^m)^-ds]\\
&+\hat{\mathbb{E}}_t[\int_0^T(n+m)|(Y_s^m)^-|^{\alpha-1}(Y_s^n)^-ds]\}].
\end{split}\end{equation}

For $2\leq \alpha<\beta$, there exist $\alpha'$, $p$, $q$, $r$, $p'$, $q'>1$, such that $\frac{1}{p}+\frac{1}{q}+\frac{1}{r}=1$, $\frac{1}{p'}+\frac{1}{q'}=1$, $(\alpha-2)\alpha' p<\beta$, $\alpha' q<\beta$, $\alpha' r<\beta$, $(\alpha-1)\alpha' p'<\beta$ and $\alpha' q'<\beta$. Applying Lemma \ref{the1.11}, Lemma \ref{the1.12} and the H\"{o}lder inequality, there exists a constant $C$ independent of $m,n$ such that
\begin{equation}\begin{split}\label{eq1.7}
&\quad \hat{\mathbb{E}}[(\int_0^T n|(Y_s^n)^-|^{\alpha-1}(Y_s^m)^-ds)^{\alpha'}]\\
&\leq \hat{\mathbb{E}}[\sup_{s\in[0,T]}\{|(Y_s^n)^-|^{(\alpha-2)\alpha'}|(Y_s^m)^-|^{\alpha'}\}(\int_0^T n(Y_s^n)^-ds)^{\alpha'}]\\
&\leq (\hat{\mathbb{E}}[\sup_{s\in[0,T]}|(Y_s^n)^-|^{(\alpha-2)\alpha' p}])^{\frac{1}{p}}(\hat{\mathbb{E}}[\sup_{s\in[0,T]}|(Y_s^m)^-|^{\alpha' q}])^{\frac{1}{q}}(\hat{\mathbb{E}}[(\int_0^T n(Y_s^n)^-ds)^{\alpha' r}])^{\frac{1}{r}}\\
&\leq C(\hat{\mathbb{E}}[\sup_{s\in[0,T]}|(Y_s^m)^-|^{\alpha' q}])^{\frac{1}{q}},
\end{split}\end{equation}
and
\begin{equation}\begin{split}\label{eq1.8}
&\quad \hat{\mathbb{E}}[(\int_0^T m|(Y_s^n)^-|^{\alpha-1}(Y_s^m)^-ds)^{\alpha'}]\\
&\leq \hat{\mathbb{E}}[\sup_{s\in[0,T]}|(Y_s^n)^-|^{(\alpha-1)\alpha'}(\int_0^T m(Y_s^m)^-ds)^{\alpha'}]\\
&\leq (\hat{\mathbb{E}}[\sup_{s\in[0,T]}|(Y_s^n)^-|^{(\alpha-1)\alpha' p'}])^{\frac{1}{p'}}(\hat{\mathbb{E}}[(\int_0^T m(Y_s^m)^-ds)^{\alpha' q'}])^{\frac{1}{q'}}\\
&\leq C(\hat{\mathbb{E}}[\sup_{s\in[0,T]}|(Y_s^n)^-|^{(\alpha-1)\alpha' p'}])^{\frac{1}{p'}}.
\end{split}\end{equation}
Then by Theorem \ref{the1.2} and Lemma \ref{the1.12},  inequalities \eqref{eq1.6}-\eqref{eq1.8} yield that
\begin{displaymath}
\lim_{n,m\rightarrow\infty}\hat{\mathbb{E}}[\sup_{t\in[0,T]}|Y_t^n-Y_t^m|^\alpha]=0.
\end{displaymath}
\end{proof}

\section{Existence and uniqueness of reflected $G$-BSDE with a lower obstacle}
\begin{theorem}\label{the1.14}
Suppose that $\xi$, $f$ satisfy (H1)-(H3) and $S$ satisfies (H4) or (H4'). Then the reflected $G$-BSDE with data $(\xi,f,S)$ has a unique solution $(Y,Z,A)$. Moreover, for any $2\leq \alpha<\beta$ we have $Y\in S^\alpha_G(0,T)$, $Z\in H_G^\alpha(0,T)$ and $A\in S_G^{\alpha}(0,T)$.
\end{theorem}

\begin{proof}
The uniqueness of the solution is a direct consequence of the a priori estimates in Proposition \ref{the1.7}, Proposition \ref{the1.8} and Proposition \ref{the1.10}.

To prove the existence, it suffices to prove the $S\equiv 0$ case. Recalling penalized $G$-BSDEs \eqref{eq1.3}, set $\hat{Y}_t=Y_t^n-Y_t^m, \hat{Z}_t=Z_t^n-Z_t^m, \hat{K}_t=K_t^n-K_t^m, \hat{L}_t=L_t^n-L_t^m$ and $\hat{f}_t=f(t,Y_t^n,Z_t^n)-f(t,Y_t^m,Z_t^m)$. By Lemma \ref{the1.13}, there exists $Y\in S_G^\alpha(0,T)$ satisfying $\lim_{n\rightarrow\infty}\hat{\mathbb{E}}[\sup_{t\in[0,t]}|Y_t-Y_t^n|^\alpha]=0$. Applying It\^{o}'s formula to $|\hat{Y}_t|^2$, we get
\begin{displaymath}
\begin{split}
&\quad|\hat{Y}_t|^2+\int_t^T|\hat{Z}_s|^2 d\langle B\rangle_s\\
&=\int_t^T2\hat{Y}_s\hat{f}_sds-\int_t^T 2\hat{Y}_sd\hat{K}_s
+\int_t^T2\hat{Y}_sd\hat{L}_s-\int_t^T 2\hat{Y}_s\hat{Z}_s dB_s\\
&\leq 2L\int_t^T[|\hat{Y}_s|^2+|\hat{Y}_s||\hat{Z}_s|]ds-\int_t^T 2\hat{Y}_sd\hat{K}_s
+\int_t^T2\hat{Y}_sd\hat{L}_s-\int_t^T 2\hat{Y}_s\hat{Z}_s dB_s.
\end{split}
\end{displaymath}
Note that for each $\varepsilon>0$,
\[2L\int_t^T |\hat{Y}_s||\hat{Z}_s|ds\leq L^2/\varepsilon\int_t^T|\hat{Y}_s|^2ds+\varepsilon\int_t^T|\hat{Z}_s|^2ds .\]
Choosing $\varepsilon<\underline{\sigma}^2$, we have
\begin{equation}\label{e1}
\begin{split}
\int_0^T|\hat{Z}_s|^2 ds
&\leq C(\int_0^T|\hat{Y}_s|^2ds-\int_0^T \hat{Y}_sd\hat{K}_s
+\int_0^T\hat{Y}_sd\hat{L}_s-\int_0^T \hat{Y}_s\hat{Z}_s dB_s)\\
&\leq C(\sup_{s\in[0,T]}|\hat{Y}_s|^2+\sup_{s\in[0,T]}|\hat{Y}_s|(|K_T^n|+|K_T^m|+|L_T^n|+|L_T^m|)-\int_0^T \hat{Y}_s\hat{Z}_s dB_s).
\end{split}
\end{equation}
By Proposition \ref{the1.3}, for any $\varepsilon'>0$, we obtain
\begin{align*}
\hat{\mathbb{E}}[(\int_0^T \hat{Y}_s\hat{Z}_sdB_s)^{\frac{\alpha}{2}}]\leq &C\hat{\mathbb{E}}[(\int_0^T \hat{Y}_s^2\hat{Z}_s^2 ds)^{\frac{\alpha}{4}}]\\
\leq&C(\hat{\mathbb{E}}[\sup_{t\in[0,T]}|\hat{Y}_t|^\alpha])^{1/2}(\hat{\mathbb{E}}[(\int_0^T |\hat{Z}_s|^2ds)^{\frac{\alpha}{2}}])^{1/2}\\
\leq &\frac{C}{4\varepsilon'}\hat{\mathbb{E}}[\sup_{t\in[0,T]}|\hat{Y}_t|^\alpha]+C\varepsilon'\hat{\mathbb{E}}[(\int_0^T |\hat{Z}_s|^2ds)^{\frac{\alpha}{2}}].
\end{align*}
Applying Lemma \ref{the1.11} and the H\"{o}lder inequality, choosing $\varepsilon'$ small enough, it follows from \eqref{e1} that
\[\hat{\mathbb{E}}[(\int_0^T|Z_s^n-Z_s^m|^2ds)^{\frac{\alpha}{2}}]\leq C\{\hat{\mathbb{E}}[\sup_{t\in[0,T]}|\hat{Y}_t|^\alpha]+(\hat{\mathbb{E}}[\sup_{t\in[0,T]}|\hat{Y}_t|^\alpha])^{1/2}\}.\]
It is straightforward to show that
\begin{displaymath}
\lim_{n,m\rightarrow \infty}\hat{\mathbb{E}}[(\int_0^T|Z_s^n-Z_s^m|^2ds)^{\frac{\alpha}{2}}]=0.
\end{displaymath}
Then there exists  a process $\{Z_t\}\in H_G^\alpha(0,T)$ such that $\hat{\mathbb{E}}[(\int_0^T|Z_s-Z_s^n|^2ds)^{\alpha/2}]\rightarrow 0$ as $n\rightarrow \infty$. Set $A_t^n=L_t^n -K_t^n$, it is easy to check that $(A_t^n)_{t\in[0,T]}$ is a nondecreasing process and
\begin{displaymath}
A_t^n-A_t^m=\hat{Y}_0-\hat{Y}_t-\int_0^t \hat{f}_sds+\int_0^t \hat{Z}_sdB_s.
\end{displaymath}
By applying Proposition \ref{the1.3} and the assumption of $f$, it follows that
\begin{displaymath}
\begin{split}
\hat{\mathbb{E}}[\sup_{t\in[0,T]}|A_t^n-A_t^m|^\alpha]
&\leq C\hat{\mathbb{E}}[\sup_{t\in[0,T]}|\hat{Y}_t|^\alpha+(\int_0^T|\hat{f}_s| ds)^\alpha
+\sup_{t\in[0,T]}|\int_0^t \hat{Z}_sdB_s|^\alpha]\\
&\leq C\{\hat{\mathbb{E}}[\sup_{t\in[0,T]}|\hat{Y}_t|^\alpha]+\hat{\mathbb{E}}[(\int_0^T|\hat{Z}_s|^2ds)^{\alpha/2}]\}\rightarrow 0.
\end{split}
\end{displaymath}
Then there exists a nondecreasing process $(A_t)_{t\in[0,T]}$ satisfying that
\begin{displaymath}
\lim_{n\rightarrow \infty}\hat{\mathbb{E}}[\sup_{t\in[0,T]}|A_t-A_t^n|^\alpha]=0.
\end{displaymath}



In the following it remains to prove that $Y_t\geq 0$, $t\in[0,T]$ and $\{-\int_0^t Y_sdA_s\}_{t\in[0,T]}$ is a non-increasing $G$-martingale. For the first statement, it can be deduced easily from Lemma \ref{the1.12}. Set
\begin{displaymath}
\widetilde{K}_t^n:=\int_0^t Y_sdK_s^n.
\end{displaymath}
Since $Y_t\geq 0$, $\forall 0\leq t\leq T$ and $K^n$ is a decreasing $G$-martingale, then $\widetilde{K}^n$ is a decreasing $G$-martingale.  Note that
\begin{align*}
\sup_{t\in[0,T]}|-\int_0^t Y_sdA_s-\widetilde{K}_t^n|
\leq &\sup_{t\in[0,T]}\{|-\int_0^t Y_sdA_s+\int_0^tY_sdA_s^n|+|\int_0^t (Y_s^n -Y_s)dA_s^n|\\
&+|\int_0^t (Y_s^n -Y_s) dK_s^n|+|\int_0^t -Y_s^n n(Y_s^n)^- ds|\}\\
\leq &\sup_{t\in[0,T]}\{|\int_0^t\widetilde{Y}_s^m d(A_s^n -A_s)|+|\int_0^t (Y_s-\widetilde{Y}_s^m)d(A_s^n -A_s)|\}\\
&+\sup_{t\in[0,T]}|Y_s-Y_s^n|[|A_T^n|+|K_T^n|]+\sup_{t\in[0,T]}(Y_s^n)^-|L_T^n|\\
=& I+II+III+IV,
\end{align*}
where $\widetilde{Y}_t^m=\sum_{i=0}^{m-1}Y_{t_i^m} I_{[t_{i}^m,t_{i+1}^m)}(t)$ and $t_i^m=\frac{iT}{m}$, $i=0,1,\cdots,m$. By simple analysis, we have
\begin{align*}
\hat{\mathbb{E}}[I]
\leq& \sum_{i=0}^{m-1}\hat{\mathbb{E}}[\sup_{s\in[0,T]}|Y_s|(|A^n_{t_{i+1}^m}-A_{t_{i+1}^m}|+|A^n_{t_{i}^m}-A_{t_{i}^m}|)]\\
\leq &(\hat{\mathbb{E}}[\sup_{s\in[0,T]}|Y_s|^2])^{1/2}
\sum_{i=0}^{m-1}\{(\hat{\mathbb{E}}[|A^n_{t_{i+1}^m}-A_{t_{i+1}^m}|^2])^{1/2}+(\hat{\mathbb{E}}[|A^n_{t_{i}^m}-A_{t_{i}^m}|^2])^{1/2}\},\\
\hat{\mathbb{E}}[II]
\leq& (\hat{\mathbb{E}}[\sup_{t\in[0,T]}|Y_s-\widetilde{Y}_s^m|^2])^{1/2}\{
(\hat{\mathbb{E}}[|A_T^n|^2])^{1/2}+(\hat{\mathbb{E}}[|A_T|^2])^{1/2}\},\\
\hat{\mathbb{E}}[III]\leq &(\hat{\mathbb{E}}[\sup_{t\in[0,T]}|Y_s-Y_s^n|^2])^{1/2}\{(\hat{\mathbb{E}}[|A_T^n|^2])^{1/2}+(\hat{\mathbb{E}}[|K_T^n|^2])^{1/2}\},\\
\hat{\mathbb{E}}[IV]\leq &(\hat{\mathbb{E}}[\sup_{t\in[0,T]}|(Y_s^n)^-|^2])^{1/2}(\hat{\mathbb{E}}[|L_T^n|^2])^{1/2}.
\end{align*}
Then for each fixed $m$, first let $n$ tends to infinity, we conclude that
\begin{displaymath}
\lim_{n\rightarrow\infty}\hat{\mathbb{E}}[\sup_{t\in[0,T]}|-\int_0^t Y_sdA_s-\widetilde{K}_t^n|]\leq C(\hat{\mathbb{E}}[\sup_{t\in[0,T]}|Y_s-\widetilde{Y}_s^m|^2])^{1/2}.
\end{displaymath}
By Lemma 3.2 in \cite{HJPS1}, sending $m$ tends to infinity, we get $\lim_{n\rightarrow\infty}\hat{\mathbb{E}}[\sup_{t\in[0,T]}|-\int_0^t Y_sdA_s-\widetilde{K}_t^n|]=0$. It follow that $\{-\int_0^t Y_sdA_s\}$ is a non-increasing $G$-martingale.
\end{proof}

Furthermore, we have the following result.
\begin{theorem}\label{the1.15}
Suppose that $\xi$, $f$ and $g$ satisfy (H1)-(H3), $S$ satisfies (H4) or (H4'). Then the reflected $G$-BSDE with data $(\xi,f,g,S)$ has a unique solution $(Y,Z,A)$. Moreover, for any $2\leq \alpha<\beta$ we have $Y\in S^\alpha_G(0,T)$, $Z\in H_G^\alpha(0,T)$ and $A\in S_G^{\alpha}(0,T)$.
\end{theorem}
\begin{proof}
The proof is similar to that of Theorem \ref{the1.14}.
\end{proof}

We next prove a comparison theorem, similar to that of \cite{HJPS2} for non-reflected $G$-BSDEs. The proof is based on the approximation method via penalization.
\begin{theorem}\label{the1.16}
Let $(\xi^1,f^1,g^1,S^1)$ and $(\xi^2,f^2,g^2,S^2)$ be two sets of data. Suppose $S^i$ satisfy (H4) or (H4'), $\xi^i$, $f^i$ and $g^i$ satisfy (H1)-(H3) for $i=1,2$. We further assume in addition the following:
\begin{description}
\item[(i)] $\xi^1\leq \xi^2$, $q.s.$;
\item[(ii)] $f^1(t,y,z)\leq f^2(t,y,z)$, $g^1(t,y,z)\leq g^2(t,y,z)$, $\forall (y,z)\in\mathbb{R}^2$;
\item[(iii)] $S_t^1\leq S^2_t$, $0\leq t\leq T$, $q.s.$.
\end{description}
Let $(Y^i,Z^i,A^i)$ be a solution of the reflected $G$-BSDE with data $(\xi^i,f^i,g^i,S^i)$, $i=1,2$ respectively. Then
\begin{displaymath}
Y_t^1\leq Y^2_t, \quad 0\leq t\leq T  \quad q.s.
\end{displaymath}
\end{theorem}

\begin{proof}
We first consider the following $G$-BSDEs parameterized by $n=1,2,\cdots$,
\begin{displaymath}
y_t^n=\xi^1+\int_t^T f^1(s,y_s^n,z_s^n)ds+\int_t^T g^1(s,y_s^n,z_s^n)d\langle B\rangle_s+\int_t^T n(y_s^n-S_s^1)^-ds-\int_t^T z_s^ndB_s-(K_T^n-K_t^n).
\end{displaymath}
By a similar analysis as the proof of Theorem \ref{the1.14}, it follows that $\lim_{n\rightarrow \infty}\hat{\mathbb{E}}[\sup_{t\in[0,T]}|Y_t^1-y_t^n|^\alpha]=0$, where $2\leq \alpha<\beta$. Noting that $(Y^2,Z^2,A^2)$ is the solution of the reflected $G$-BSDE with data $(\xi^2,f^2,g^2,S^2)$ and $Y^2_t\geq S^2_t$, $0\leq t\leq T$, we have
\begin{displaymath}
Y^2_t=\xi^2+\int_t^T f^2(s,Y^2_s,Z^2_s)ds+\int_t^T g^2(s,Y_s^2,Z_s^2)d\langle B\rangle_s+\int_t^T n(Y^2_s-S^2_s)^-ds-\int_t^T Z^2_sdB_s+(A^2_T-A^2_t).
\end{displaymath}
Applying Theorem \ref{the1.5} yields $Y^2_t\geq y_t^n$, for all $n\in\mathbb{N}$. Letting $n\rightarrow \infty$, we conclude that $Y^2_t\geq Y^1_t$.
\end{proof}

\section{Relation between  reflected $G$-BSDEs and obstacle problems for nonlinear parabolic PDEs}

In this section, we will give a probabilistic representation for solutions of some obstacle problems for fully nonlinear parabolic PDEs using the reflected $G$-BSDE we have mentioned in the above sections. For this purpose, we need to put the reflected $G$-BSDE in a Markovian framework.

For each $(t,x)\in[0,T]\times\mathbb{R}^d$, let $\{X_s^{t,x}, t\leq s\leq T\}$ be the unique $\mathbb{R}^d$-valued solution of the SDE driven by $G$-Brownian motion (here we use Einstein convention):
\begin{equation}\label{eq1.6}
X_s^{t,x}=x+\int_t^s b(r,X_r^{t,x})dr+\int_t^s l_{ij}(r,X_r^{t,x})d\langle B^i,B^j\rangle_r+\int_t^s \sigma_{i}(r,X_r^{t,x})dB_r^i.
\end{equation}
We assume that the data $(\xi,f,g,S)$ of the reflected $G$-BSDE take the following form:
\begin{align*}
\xi&=\phi(X_T^{t,x}),\ \ \ \ \ \  f(s,y,z)=f(s,X_s^{t,x},y,z),\\
S_s&=h(s,X_s^{t,x}), \ \ g_{ij}(s,y,z)=g_{ij}(s,X_s^{t,x},y,z),
\end{align*}
where $b:[0,T]\times\mathbb{R}^d\rightarrow \mathbb{R}^d$, $l_{ij}:[0,T]\times\mathbb{R}^d\rightarrow \mathbb{R}^{d}$, $\sigma_i:[0,T]\times\mathbb{R}^d\rightarrow \mathbb{R}^{d}$, $\phi:\mathbb{R}^d\rightarrow \mathbb{R}$, $f,g_{ij}:[0,T]\times\mathbb{R}^d\times\mathbb{R}\times\mathbb{R}^d\rightarrow \mathbb{R}$ and $h:[0,T]\times\mathbb{R}^d\rightarrow \mathbb{R}$ are deterministic functions and satisfy the following conditions:
\begin{description}
\item[(A1)] $l_{ij}=l_{ji}$ and $g_{ij}=g_{ji}$ for $1\leq i,j\leq d$;
\item[(A2)] $b$, $l_{ij}$, $\sigma_i$, $f$, $g_{ij}$ are continuous in $t$;
\item[(A3)] There exist a positive integer $m$ and a constant $L$ such that
\begin{align*}
&|b(t,x)-b(t,x')|+\sum_{i,j=1}^d|l_{ij}(t,x)-l_{ij}(t,x')|+\sum_{i=1}^d|\sigma_i(t,x)-\sigma_i(t,x')|\leq L|x-x'|,\\
&|\phi(x)-\phi(x')|\leq L(1+|x|^m+|x'|^m)|x-x'|,\\
&|f(t,x,y,z)-f(t,x',y',z')|+\sum_{i.j=1}^d|g_{ij}(t,x,y,z)-g_{ij}(t,x',y',z')|\\
&\leq L[(1+|x|^m+|x'|^m)|x-x'|+|y-y'|+|z-z'|].
\end{align*}
\item[(A4)]  $h$ is uniformly continuous w.r.t $(t,x)$ and bounded from above, $h(T,x)\leq \Phi(x)$ for any $x\in\mathbb{R}^d$;
\item[(A4$^{\prime}$)]$h$ belongs to the space $C^{1,2}_{Lip}([0,T]\times \mathbb{R}^d)$ and $h(T,x)\leq \Phi(x)$ for any $x\in\mathbb{R}^d$, where
$C^{1,2}_{Lip}([0,T]\times \mathbb{R}^d)$ is the space of  all functions of class $C^{1,2}([0,T]\times\mathbb{R}^d)$ whose partial derivatives of order less than or equal to $2$  and itself are  Lipschtiz continuous functions with respect to $x$.
\end{description}

We have the following estimates of $G$-SDEs, which come from Chapter V of Peng \cite{P10}.
\begin{proposition}[\cite{P10}]\label{the1.17}
Let $\xi,\xi'\in L_G^p(\Omega_t;\mathbb{R}^d)$ and $p\geq 2$. Then we have, for each $\delta\in[0,T-t]$,
\begin{align*}
\hat{\mathbb{E}}_t[\sup_{s\in[t,t+\delta]}|X_{s}^{t,\xi}-X_{s}^{t,\xi'}|^p]&\leq C|\xi-\xi'|^p,\\
\hat{\mathbb{E}}_t[|X_{t+\delta}^{t,\xi}|^p]&\leq C(1+|\xi|^p),\\
\hat{\mathbb{E}}_t[\sup_{s\in[t,t+\delta]}|X_s^{t,\xi}-\xi|^p]&\leq C(1+|\xi|^p)\delta^{p/2},
\end{align*}
where the constant $C$ depends on $L,G,p,d$ and $T$.
\end{proposition}

\begin{proof}
For the reader's convenience, we give a brief proof here. It is easy to check that $\{X_s^{t,\xi}\}_{s\in[t,T]}$, $\{X_s^{t,\xi'}\}_{s\in[t,T]}\in M_G^p(0,T;\mathbb{R}^d)$. By Proposition \ref{the1.3}, we have
\begin{align*}
&\hat{\mathbb{E}}_t[\sup_{s\in[t,t+\delta]}| X_s^{t,\xi}-X_s^{t,\xi'}|^p]\\
\leq &C\hat{\mathbb{E}}_t[|\xi-\xi'|^p+\int_t^{t+\delta}|X_s^{t,\xi}-X_s^{t,\xi'}|^pds+\sup_{s\in[t,t+\delta]}|\int_t^s (\sigma(r,X_r^{t,\xi})-\sigma(r,X_r^{t,\xi'}))dB_r|^p]\\
\leq& C\{|\xi-\xi'|^p+\hat{\mathbb{E}}_t[\int_t^{t+\delta}|X_s^{t,\xi}-X_s^{t,\xi'}|^pds]+\hat{\mathbb{E}}_t[(\int_t^{t+\delta}|X_s^{t,\xi}-X_s^{t,\xi'}|^2ds)^{p/2}]\}\\
\leq &C\{|\xi-\xi'|^p+\int_t^{t+\delta}\hat{\mathbb{E}}_t[\sup_{r\in[t,s]}|X_r^{t,\xi}-X_r^{t,\xi'}|^p]ds\}.
\end{align*}
By the Gronwall inequality, we get the first inequality. The others can be proved similarly.
\end{proof}

It follows from the previous results that for each $(t,x)\in[0,T]\times\mathbb{R}^d$, there exists a unique triple $(Y^{t,x}_s,Z^{t,x}_s,A^{t,x}_s)_{s\in[t,T]}$, which solves the following reflected $G$-BSDE:
\begin{description}
\item[(i)] $Y_s^{t,x}=\phi(X_T^{t,x})+\int_s^T f(r,X_r^{t,x},Y_r^{t,x},Z_r^{t,x})dr+\int_s^T g_{ij}(r,X_r^{t,x},Y_r^{t,x},Z_r^{t,x})d\langle B^i,B^j\rangle_r\\
\ \ \ \ \ \ -\int_s^T Z_r^{t,x}dB_r+A_T^{t,x}-A_s^{t,x},\quad t\leq s\leq T$;
\item[(ii)] $Y_s^{t,x}\geq h(s,X_s^{t,x})$, $t\leq s\leq T$;
\item[(iii)] $\{A_s^{t,x}\}$ is nondecreasing and continuous, and $\{-\int_t^s (Y_r^{t,x}-h(r,X_r^{t,x}))dA_r^{t,x}, t\leq s\leq T\}$ is a non-increasing $G$-martingale.
\end{description}

We now consider the following obstacle problem for a parabolic PDE.
\begin{equation}\label{eq1.7}
\begin{cases}
\min(-\partial_t u(t,x)-F(D_x^2 u,D_x u,u,x,t),u(t,x)-h(t,x))=0,
 &(t,x)\in(0,T)\times \mathbb{R}^d,\\
u(T,x)=\phi(x),  &x\in\mathbb{R}^d,
\end{cases}
\end{equation}
where
\begin{align*}
F(D_x^2 u,D_x u,u,x,t)=&G(H(D_x^2 u,D_x u,u,x,t))+\langle b(t,x),D_x u\rangle\\
&+f(t,x,u,\langle \sigma_1(t,x),D_x u\rangle,\cdots,\langle \sigma_d(t,x),D_x u\rangle),\\
H(D_x^2 u,D_x u,u,x,t)=&\langle D_x^2u\sigma_i(t,x),\sigma_j(t,x)\rangle+2\langle D_xu,l_{ij}(t,x)\rangle\\
&+2g_{ij}(t,x,u,\langle \sigma_1(t,x),D_x u\rangle,\cdots,\langle \sigma_d(t,x),D_x u\rangle).
\end{align*}

We need to consider solutions of the above PDE in the viscosity sense. The best candidate to define the notion of viscosity solution is by using the language of sub- and super-jets; (see \cite{CIL}).

\begin{definition}
Let $u\in C((0,T)\times\mathbb{R}^d)$ and $(t,x)\in(0,T)\times \mathbb{R}^d$. We denote by $\mathcal{P}^{2,+} u(t,x)$ [the ``parabolic superjet" of $u$ at $(t,x)$] the set of triples $(p,q,X)\in\mathbb{R}\times\mathbb{R}^d\times \mathbb{S}_d$ which are such that
\begin{align*}
u(s,y)\leq& u(t,x)+p(s-t)+\langle q,y-x\rangle\\
&+\frac{1}{2}\langle X(y-x),y-x\rangle+o(|s-t|+|y-x|^2).
\end{align*}
Similarly, we define $\mathcal{P}^{2,-} u(t,x)$ [the ``parabolic subjet" of $u$ at $(t,x)$] by $\mathcal{P}^{2,-} u(t,x):=-\mathcal{P}^{2,+}(- u)(t,x)$.
\end{definition}

Then we can give the definition of the viscosity solution of the obstacle problem \eqref{eq1.7}.

\begin{definition}
It can be said that $u\in C([0,T]\times\mathbb{R}^d)$ is a viscosity subsolution of \eqref{eq1.7} if $u(T,x)\leq \phi(x)$, $x\in\mathbb{R}^d$, and at any point $(t,x)\in(0,T)\times\mathbb{R}^d$, for any $(p,q,X)\in\mathcal{P}^{2,+}u(t,x)$,
\begin{displaymath}
\min(u(t,x)-h(t,x), -p-F(X,q,u(t,x),x,t))\leq0.
\end{displaymath}
It can be said that $u\in C([0,T]\times\mathbb{R}^d)$ is a viscosity supersolution of \eqref{eq1.7} if $u(T,x)\geq\phi(x)$, $x\in\mathbb{R}^d$, and at any point $(t,x)\in(0,T)\times\mathbb{R}^d$, for any $(p,q,X)\in\mathcal{P}^{2,-}u(t,x)$,
\begin{displaymath}
\min(u(t,x)-h(t,x), -p-F(X,q,u(t,x),x,t))\geq 0.
\end{displaymath}
$u\in C([0,T]\times\mathbb{R}^d)$ is said to be a viscosity solution of \eqref{eq1.7} if it is both a viscosity subsolution and supersolution.
\end{definition}

We now define
\begin{equation}\label{eq1.8}
u(t,x):=Y_t^{t,x},\quad (t,x)\in[0,T]\times\mathbb{R}^d.
\end{equation}
It is important to note that $u(t,x)$ is a deterministic function. We claim that $u$ is a continuous function. For simplicity, we only consider the case $g=0$ in the next three lemmas. The results still hold for the other cases. 

\begin{lemma}\label{the1.18}
Let assumption (A1)-(A3) and (A4') hold. For each $t\in[0,T]$, $x_1,x_2\in\mathbb{R}^d$, we have
\begin{displaymath}
|u(t,x_1)-u(t,x_2)|\leq C(1+|x_1|^{m\vee2}+|x_2|^{m\vee2})|x_1-x_2|.
\end{displaymath}
\end{lemma}

\begin{proof}
From Proposition \ref{the1.8}, since $u(t,x)$ is a deterministic function, we have
\begin{equation}
\begin{split}
|u(t,x_1)-u(t,x_2)|^2&\leq C\{\hat{\mathbb{E}}[|(\phi(X_T^{t,x_1})-h(T,X_T^{t,x_1}))-(\phi(X_T^{t,x_2})-h(T,X_T^{t,x_2}))|^2\\
&\quad+\int_t^T|f(s,X_s^{t,x_1},Y_s^{t,x_1},Z_s^{t,x_1})-f(s,X_s^{t,x_2},Y_s^{t,x_1},Z_s^{t,x_1})|^2 ds\\
&\quad +\int_t^T|b^1(s)-b^2(s)|^2+|l^1_{ij}(s)-l^2_{ij}(s)|^2+|\sigma^1_i(s)-\sigma^2_i(s)|^2\\
&\quad+|h(s,X_s^{t,x_1})-h(s,X_s^{t,x_2})|^2 ds]+|h(t,x_1)-h(t,x_2)|^2\},
\end{split}
\end{equation}
where for $k=1,2$,
\begin{align*}
&{b}^k(s)=\partial_s h(s,X_s^{t,x_k})+\langle b(s,X_s^{t,x_k}),D_x h(s,X_s^{t,x_k})\rangle,\\
&{l}^k_{ij}(s)=\langle D_x h(s,X_s^{t,x_k}), l_{ij}(s,X_s^{t,x_k})\rangle+\frac{1}{2}\langle D_x^2 h(s,X_s^{t,x_k})\sigma_i(s,X_s^{t,x_k}),\sigma_j(s,X_s^{t,x_k})\rangle,\\
&{\sigma}^k_i(s)=\langle \sigma_i(s,X_s^{t,x_k}),D_x h(s,X_s^{t,x_k})\rangle.
\end{align*}
Set $\hat{X}_s^t=X_s^{t,x_1}-X_s^{t,x_2}$. By the assumptions (A3), (A4') and Proposition \ref{the1.17}, we have
\begin{align*}
|u(t,x_1)-u(t,x_2)|^2
\leq &C\{\hat{\mathbb{E}}[(1+\sum_{k=1}^2|X_T^{t,x_k}|^m)^2|\hat{X}_T^t|^2]
+\int_t^T\hat{\mathbb{E}}[(1+\sum_{k=1}^2|X_s^{t,x_k}|^m)^2|\hat{X}_s^t|^2]ds\\
&+\int_t^T\hat{\mathbb{E}}[(1+\sum_{k=1}^2|X_s^{t,x_k}|^2)^2|\hat{X}_s^t|^2]ds
+\int_t^T\hat{\mathbb{E}}[|\hat{X}_s^t|^2]ds+|x_1-x_2|^2\}\\
\leq &C(1+|x_1|^{2m\vee4}+|x_2|^{2m\vee4})\{(\hat{\mathbb{E}}[\sup_{s\in[t,T]}|\hat{X}_s^t|^4])^{1/2}+|x_1-x_2|^2\}\\
\leq &C(1+|x_1|^{2m\vee4}+|x_2|^{2m\vee4})|x_1-x_2|^2.
\end{align*}
The proof is complete.
\end{proof}

\begin{lemma}\label{the1.19}
Let assumption (A1)-(A4) hold. For each $t\in[0,T]$, $x,x'\in\mathbb{R}^d$, we have
\begin{displaymath}
|u(t,x_1)-u(t,x_2)|^2\leq C\{(1+|x_1|^{2m}+|x_2|^{2m})|x_1-x_2|^2+(1+|x_1|^{m+1}+|x_2|^{m+1})|x_1-x_2|\}.
\end{displaymath}
\end{lemma}

\begin{proof}
From Proposition \ref{the1.10} and Proposition \ref{the1.17}, by a similar analysis with the above lemma, we get the desired result.
\end{proof}

The following lemma states that $u(t,x)$ is continuous with respect to $t$.
\begin{lemma}\label{the1.20}
The function $u(t,x)$ is continuous in $t$.
\end{lemma}

\begin{proof}
We only need to prove the case where (A1)-(A3) and (A4') hold. The case that (A1)-(A4) hold can be proved in a similar way. We define $X_s^{t,x}:=x$, $Y_s^{t,x}:=Y_t^{t,x}$, $Z_s^{t,x}:= 0$ and $A_s^{t,x}:=0$ for $0\leq s\leq t$. Then we define the obstacle
 \begin{displaymath}
 \tilde{S}^{t,x}_u=
 \begin{cases}h(t,x)+\int_t^u \tilde{b}(s,X_s^{t,x})ds+\int_t^u \tilde{l}_{ij}(s,X_s^{t,x})d\langle B^i ,B^j\rangle_s+\int_t^u \tilde{\sigma}_i(s,X_s^{t,x})dB^i_s, & u\in(t,T];\\
   h(t,x),&u\in[0,t],
    \end{cases}
    \end{displaymath}
    where 
\begin{align*}
&\tilde{b}(s,X_s^{t,x})=\partial_s h(s,X_s^{t,x})+\langle b(s,X_s^{t,x}),D_x h(s,X_s^{t,x})\rangle,\\
&\tilde{l}_{ij}(s,X_s^{t,x})=\langle D_x h(s,X_s^{t,x}), l_{ij}(s,X_s^{t,x})\rangle+\frac{1}{2}\langle D_x^2 h(s,X_s^{t,x})\sigma_i(s,X_s^{t,x}),\sigma_j(s,X_s^{t,x})\rangle,\\
&\tilde{\sigma}_i(s,X_s^{t,x})=\langle \sigma_i(s,X_s^{t,x}),D_x h(s,X_s^{t,x})\rangle.
\end{align*}
It is easy to check that $(Y^{t,x}_s,Z^{t,x}_s,A^{t,x}_s)_{s\in[0,T]}$ is the solution to the reflected $G$-BSDE with data $(\phi(X_T^{t,x}),\tilde{f}^{t,x},\tilde{S}^{t,x})$ where
\begin{align*}
\tilde{f}^{t,x}(s,y,z):=I_{[t,T]}(s)f(s,X_s^{t,x},y,z).
\end{align*}
Fix $x\in\mathbb{R}^d$, for $0\leq t_1\leq t_2\leq T$, by Proposition \ref{the1.8}, we have
\begin{align*}
&|u(t_1,x)-u(t_2,x)|^2=|Y_0^{t_1,x}-Y_0^{t_2,x}|^2\\
\leq &C\{\hat{\mathbb{E}}[|(\phi(X_T^{t_1,x})-h(T,X_T^{t_1,x}))-(\phi(X_T^{t_2,x})-h(T,X_T^{t_2,x}))|^2+|h(t_1,x)-h(t_2,x)|^2\\
&+\int_0^T|\hat{\lambda}_{t_1,t_2}(s)|^2+|\hat{\rho}_{t_1,t_2}(s)|^2+|h(s,X_s^{t_1,x})-h(s,X_s^{t_2,x})|^2ds]\},
\end{align*}
where \begin{align*} \hat{\lambda}_{t_1,t_2}(s)=|I_{[t_1,T]}(s)f(s,X_s^{t_1,x},Y_s^{t_2,x},Z_s^{t_2,x})-I_{[t_2,T]}(s)f(s,X_s^{t_2,x},Y_s^{t_2,x},Z_s^{t_2,x})|,
\end{align*}
and
\begin{align*}\hat{\rho}_{t_1,t_2}(s)=&|I_{[t_1,T]}(s)\tilde{b}(s,X_s^{t_1,x})-I_{[t_2,T]}(s)\tilde{b}(s,X_s^{t_2,x})|\\
&+|I_{[t_1,T]}(s)\tilde{l}_{ij}(s,X_s^{t_1,x})-I_{[t_2,T]}(s)\tilde{l}_{ij}(s,X_s^{t_2,x})|\\
&+|I_{[t_1,T]}(s)\tilde{\sigma}_i(s,X_s^{t_1,x})-I_{[t_2,T]}(s)\tilde{\sigma}_i(s,X_s^{t_2,x})|.
\end{align*}
Set $\hat{X}_s^x=X^{t_1,x}_s-X^{t_2,x}_s$. By H\"{o}lder's inequality  and assumptions (A3), (A4'),  we deduce that
\begin{align*}
|u(t_1,x)-u(t_2,x)|^2
\leq &C\{\hat{\mathbb{E}}[(1+|X_T^{t_1,x}|^m+|X_T^{t_2,x}|^m)^2|\hat{X}_T^x|^2]+|h(t_1,x)-h(t_2,x)|^2\\
&+\int_{t_2}^T\hat{\mathbb{E}}[(1+|X_s^{t_1,x}|^{m\vee2}+|X_s^{t_2,x}|^{m\vee2})^2|\hat{X}_s^x|^2]ds\\
&+\int_{t_1}^{t_2}\hat{\mathbb{E}}[1+|X_s^{t_1,x}|^{(2m+2)\vee 6}+|Y_s^{t_2,x}|^2]ds\}.
\end{align*}
Note that $\hat{X}_s^x=X_s^{t_2, X_{t_2}^{t_1,x}}-X_s^{t_2,x}$, for $s\in[t_2,T]$. Applying Proposition \ref{the1.17}, it follows that
\[|u(t_1,x)-u(t_2,x)|\leq C\{(1+|x|^{(m+1)\vee 3})|t_2-t_1|^{\frac{1}{2}}+|h(t_2,x)-h(t_1,x)|\}.\]
The proof is complete.
\end{proof}

We will use the approximation of the reflected $G$-BSDE  by penalization. For each $(t,x)\in[0,T]\times\mathbb{R}^d$, $n\in\mathbb{N}$, let $\{(Y_s^{n,t,x},Z_s^{n,t,x},K_s^{n,t,x}), t\leq s
\leq T\}$ denote the solution of the $G$-BSDE:
\begin{align*}
Y_s^{n,t,x}=& \phi(X_T^{t,x})+\int_s^T f(r,X_r^{t,x},Y_r^{n,t,x},Z_r^{n,t,x})dr+\int_s^T g_{ij}(r,X_r^{t,x},Y_r^{n,t,x},Z_r^{n,t,x})d\langle B^i,B^j\rangle_r\\
&+n\int_s^T (Y_r^{n,t,x}-h(r,X_r^{t,x}))^-dr-\int_s^T Z_r^{n,t,x}dB_r-(K_T^{n,t,x}-K_s^{n,t,x}), \quad t\leq s\leq T.
\end{align*}
We define
\begin{displaymath}
u_n(t,x):=Y_t^{n,t,x}, \quad 0\leq t\leq T, x\in\mathbb{R}^d.
\end{displaymath}
By Theorem 4.5 in \cite{HJPS2}, $u_n$ is the viscosity solution of the parabolic PDE
\begin{equation}\label{eq1.9}
\begin{cases}
-\partial_t u_n(t,x)-F_n(D_x^2 u_n(t,x),D_x u_n(t,x),u_n(t,x),x,t)=0, & (t,x)\in[0,T]\times\mathbb{R}^d\\
u_n(T,x)=\phi(x), &  x\in\mathbb{R}^d,
\end{cases}\end{equation}
where
\begin{displaymath}
F_n(D_x^2 u,D_x u,u,x,t)=F(D_x^2 u,D_x u,u,x,t)+n(u-h(t,x))^-.
\end{displaymath}

\begin{theorem}\label{the1.21}
The function $u$ defined by \eqref{eq1.8} is the unique viscosity solution of the obstacle problem \eqref{eq1.7}.
\end{theorem}

\begin{proof}
From the results of the previous section, for each $(t,x)\in[0,T]\times\mathbb{R}^d$, we obtain
\begin{displaymath}
u_n(t,x)\uparrow u(t,x), \textrm{ as }n\rightarrow \infty.
\end{displaymath}
By Proposition 4.2, Theorem 4.5 in \cite{HJPS2} and Lemma \ref{the1.18}-Lemma \ref{the1.20}, $u_n$ and $u$ are continuous. Then by applying Dini's theorem, the sequence $u^n$ uniformly converges to $u$ on compact sets.

We first show that $u$ is a subsolution of \eqref{eq1.7}. For each fixed $(t,x)\in(0,T)\times\mathbb{R}^d$, let $(p,q,X)\in\mathcal{P}^{2,+} u(t,x)$. Without loss of generality, we may assume that $u(t,x)>h(t,x)$. By Lemma 6.1 in \cite{CIL}, there exists sequences
\begin{displaymath}
n_j\rightarrow \infty,\ (t_j,x_j)\rightarrow (t,x),\ (p_j,q_j,X_j)\in \mathcal{P}^{2,+} u_{n_j}(t_j,x_j),
\end{displaymath}
such that $(p_j,q_j,X_j)\rightarrow(p,q,X)$. Since $u^n$ is the viscosity solution to equation \eqref{eq1.9}, it follows that for any $j$,
\begin{displaymath}
-p_j-F_{n_j}(X_j,q_j,u_{n_j}(t_j,x_j),x_j,t_j)\leq 0.
\end{displaymath}
Note that $u_n$ is uniformly convergence on compact sets, by the assumption that $u(t,x)>h(t,x)$, we derive that for $j$ large enough $u_{n_j}(t_j,x_j)>h(t_j,x_j)$; therefore, sending $j$ goes to infinity in the above inequality yields:
\begin{displaymath}
-p-F(X,q,u(t,x),x,t)\leq 0.
\end{displaymath}
Then we conclude that $u$ is a subsolution of \eqref{eq1.7}.

It remains to prove that $u$ is a supersolution of \eqref{eq1.7}. For each fixed $(t,x)\in(0,T)\times\mathbb{R}^d$, and $(p,q,X)\in\mathcal{P}^{2,-} u(t,x)$. Noting that $\{Y_s^{t,x}\}_{s\in[t,T]}$ is the solution of reflected $G$-BSDE with data $(\xi,f,g,S)$, where $S_s=h(s,X_s^{t,x})$, we have
\[u(t,x)=Y_t^{t,x}\geq h(t,x).\]
Applying Lemma 6.1 in \cite{CIL} again, there exists sequences
\begin{displaymath}
n_j\rightarrow \infty,\ (t_j,x_j)\rightarrow (t,x),\ (p_j,q_j,X_j)\in \mathcal{P}^{2,-} u_{n_j}(t_j,x_j),
\end{displaymath}
such that $(p_j,q_j,X_j)\rightarrow(p,q,X)$. Since $u^n$ is the viscosity solution to equation \eqref{eq1.9}, we derive that for any $j$,
\begin{displaymath}
-p_j-F_{n_j}(X_j,q_j,u_{n_j}(t_j,x_j),x_j,t_j)\geq 0.
\end{displaymath}
Therefore
\begin{displaymath}
-p_j-F(X_j,q_j,u_{n_j}(t_j,x_j),x_j,t_j)\geq 0.
\end{displaymath}
Sending $j\rightarrow \infty$ in the above inequality, we have
\begin{displaymath}
-p-F(X,q,u(t,x),x,t)\geq 0,
\end{displaymath}
which implies that $u$ is a supersolution of \eqref{eq1.7}. Thus $u$ is a viscosity solution of \eqref{eq1.7}.

Analysis similar to the proof of Theorem 8.6 in \cite{KKPPQ} shows that there exists at most one solution of the obstacle problem \eqref{eq1.7} in the class of continuous functions which grow at most polynomially at infinity. The proof is complete.
\end{proof}

\section{American options under volatility uncertainty}
Now let us consider the financial market with volatility uncertainty. The market model $\mathcal{M}$ is introduced in \cite{Vorbrink} consisting of two assets whose dynamics are given by
\begin{align*}
&d\gamma_t=r \gamma_t dt, \quad \gamma_0=1,\\
&dS_t=rS_tdt+S_tdB_t, \quad S_0=x_0>0,
\end{align*}
where $r\geq 0$ is a constant interest rate. The asset $\gamma=(\gamma_t)$ represents a riskless bond. The stock price is described by a geometric $G$-Brownian motion.  Since the deviation of the process $B$ from its mean is unknown, this model shows the ambiguity under volatility uncertainty.

\begin{definition}[\cite{Vorbrink}]\label{def1.1}
A cumulative consumption process $C=(C_t)$ is a nonnegative $\mathcal{F}_t$-adapted process with values in $L_G^1(\Omega_T)$, and with nondecreasing, RCLL paths on $(0,T]$, and $C_0=0$, $C_T<\infty$, $q.s.$. A portfolio process $\pi=(\pi_t)$ is an $\mathcal{F}_t$-adapted real valued process with values in $L_G^1(\Omega_T)$.
\end{definition}

\begin{remark}
Let $(\Omega_T,\mathcal{F},P)$ be the classical probability space, where $\Omega_T=C_0([0,T],\mathbb{R})$ and $\mathcal{F}=\mathcal{B}(\Omega_T)$. Furthermore, there exists a process $W=\{W_t\}$ which is a Brownian motion with respect to $P$. Then the filtration in the above definition is given by $\mathcal{F}_t=\sigma\{W_s|0\leq s\leq t\}\vee \mathcal{N}$, where $\mathcal{N}$ denotes the collection of all $P$-null subsets.
\end{remark}

\begin{definition}[\cite{Vorbrink}]\label{def1.2}
For a given initial capital $y$, a portfolio process $\pi$, and a cumulative consumption process $C$, consider the wealth equation
\begin{align*}
dX_t&=X_t(1-\pi_t)\frac{d\gamma_t}{\gamma_t}+X_t\pi_t\frac{dS_t}{S_t}-dC_t\\
&=rX_tdt+\pi_tX_tdB_t-dC_t
\end{align*}
with initial wealth $X_0=y$. Or equivalently
\begin{displaymath}
\gamma_t^{-1}X_t=y-\int_0^t \gamma_u^{-1}dC_u+\int_0^t\gamma_u^{-1}X_u\pi_udB_u, \quad \forall t\leq T.
\end{displaymath}

If this equation has a unique solution $X=(X_t):=X^{y,\pi,C}$, it is called the wealth process corresponding to the triple $(y,\pi,C)$.
\end{definition}

\begin{definition}[\cite{Vorbrink}]
A portfolio/consumption process pair $(\pi,C)$ is called admissible for an initial capital $y\in\mathbb{R}$ if
\begin{description}
\item[(i)] the pair obeys the conditions of Definition \ref{def1.1} and \ref{def1.2};
\item[(ii)] $(\pi_tX_t^{y,\pi,C})\in M_G^2(0,T)$;
\item[(iii)] the solution $X_t^{y,\pi,C}$ satisfies
\begin{displaymath}
X_t^{y,\pi,C}\geq -L, \quad \forall t\leq T, q.s.
\end{displaymath}
where $L$ is a nonnegative random variable in $L_G^2(\Omega_T)$.
\end{description}
We then write $(\pi,C)\in \mathcal{A}(y)$.
\end{definition}

We denote by $\mathcal{T}_{s,t}$ the set of all stopping times taking values in $[s,t]$, for any $0\leq s\leq t\leq T$. Then the American contingent claims may be defined by the following:
\begin{definition}[\cite{Karatzas}]
An American contingent claim is a financial instrument consisting of
\begin{description}
\item[(i)] an expiration date $T\in(0,\infty]$;
\item[(ii)] the selection of an exercise time $\tau\in \mathcal{T}_{0,T}$;
\item[(iii)] a payoff $H_{\tau}$ at the exercise time.
\end{description}
\end{definition}

We should require that the payoff process $\{H_t\}_{t\in[0,T]}$ satisfy (H4) or (H4') in Section 3. Since the financial market under volatility uncertainty is incomplete, it is natural to consider the superhedging price for the American contingent claims.
\begin{definition}
Given an American contingent claim $(T,H)$ we define the superhedging class
\begin{displaymath}
\mathcal{U}:=\{y\geq 0|\exists (\pi,C)\in \mathcal{A}(y):\textrm{ for any stopping time }\tau, X_\tau^{y,\pi,C}\geq H_\tau, q.s.\}.
\end{displaymath}
The superhedging price is defined as
\begin{displaymath}
h_{up}:=\inf\{u|y\in\mathcal{U}\}.
\end{displaymath}
\end{definition}

\begin{theorem}
Given the financial market $\mathcal{M}$ and an American contingent claim $(T,H)$, we have $h_{up}=Y_0$, where $Y=(Y_t)$ is the solution to the reflected $G$-BSDE with parameter $(H_T, f, H_t)$ with $f(y)=ry$.
\end{theorem}

\begin{proof}
Let $y\in\mathcal{U}$. By the definition of $\mathcal{U}$, there exists a pair $(\pi,C)\in\mathcal{A}(y)$ such that for any stopping time $\tau$, $X_\tau^{y,\pi,C}\geq H_\tau$. Applying Lemma 3.4, Lemma 4.2 and Lemma 4.3 in \cite{LP}, we derive that for any $\eta\in M_G^2(0,T)$
\[\hat{\mathbb{E}}[\int_0^\tau \eta_s dB_s]=0.\]
Then we obtain
\begin{align*}
y&= \hat{\mathbb{E}}[y+\int_0^\tau \gamma_u^{-1}X_u^{y,\pi,C}\pi_u dB_u]\\
&\geq \hat{\mathbb{E}}[y+\int_0^\tau \gamma_u^{-1}X_u^{y,\pi,C}\pi_u dB_u-\int_0^\tau \gamma_u^{-1}dC_u]\\
&=\hat{\mathbb{E}}[\gamma_{\tau}^{-1}X_\tau^{y,\pi,C}]\geq \hat{\mathbb{E}}[\gamma_\tau^{-1}H_\tau].
\end{align*}
It follows that $h_{up}\geq \sup_{\tau\in\mathcal{T}_{0,T}}\hat{\mathbb{E}}[\gamma_\tau^{-1}H_\tau]$.

Now we turn to prove the inverse inequality. Consider the following reflected $G$-BSDE
\begin{align*}
\begin{cases}Y_t=H_T-\int_t^T rY_sds-\int_t^T Z_sdB_s+(A_T-A_t),\\
Y_t\geq H_t.
\end{cases}
\end{align*}
By Theorem \ref{the1.14}, there exists a unique solution $(Y,Z,A)$ to the above equation. Let $C=A$, $\pi=\frac{Z}{Y}$. Then $H_\tau \leq Y_\tau=X^{Y_0,\pi,C}_\tau$, which implies $Y_0\in\mathcal{U}$. It follows that $h_{up}\leq Y_0$. Applying It\^{o}'s formula to $\widetilde{Y}_t=\gamma_t^{-1}Y_t$, we conclude that $\widetilde{Y}$ is a solution to the reflected $G$-BSDE with data $(\gamma_T^{-1}H_T,0,\gamma^{-1}H)$. By the following proposition, we finally get the desired result.
\end{proof}

\begin{proposition}
Let $(Y,Z,A)$ be a solution of the reflected $G$-BSDE with data $(\xi,f,S)$. Then we have
\begin{displaymath}
Y_0=\sup_{\tau\in\mathcal{T}_{0,T}} \hat{\mathbb{E}}[\int_0^\tau f(s,Y_s,Z_s)ds+S_\tau I_{\{\tau<T\}}+\xi I_{\{\tau=T\}}].
\end{displaymath}
\end{proposition}

\begin{proof}
Let $\tau\in\mathcal{T}_{0,T}$. Note the fact that
\[\hat{\mathbb{E}}[\int_0^\tau Z_s dB_s]=0.\]
Then we have
\begin{align*}
Y_0&=\hat{\mathbb{E}}[\int_0^\tau f(s,Y_s,Z_s)ds+Y_\tau+A_\tau]\\
&\geq \hat{\mathbb{E}}[\int_0^\tau f(s,Y_s,Z_s)ds+S_\tau I_{\{\tau<T\}}+\xi I_{\{\tau=T\}}].
\end{align*}

We are now in a position to show the inverse inequality. By the definition of the solution of the reflected $G$-BSDE, we may define
\begin{displaymath}
K_t:=-\int_0^t (Y_s-S_s)dA_s,
\end{displaymath}
then $K$ is a non-increasing $G$-martingale. Let
\begin{displaymath}
D^n=\inf\{0\leq t\leq T:Y_t-S_t<\frac{1}{n}\}\wedge T.
\end{displaymath}
By example \ref{th13}, $D^n$ is a $*$-stopping time for $n\geq 1$. It is easy to check that $D^n\rightarrow D$, where
\begin{displaymath}
D=\inf\{0\leq t\leq T:Y_t-S_t=0\}\wedge T.
\end{displaymath}
Noting that $A$ is nondecreasing, by Theorem \ref{th14}, it follows that
\begin{displaymath}
0=\hat{\mathbb{E}}[K_{D^n}]=\hat{\mathbb{E}}[-\int_0^{D^n} (Y_s-S_s)dA_s]\leq \frac{1}{n}\hat{\mathbb{E}}[-A_{D^n}]\leq  0,
\end{displaymath}
which yields $\hat{\mathbb{E}}[-A_{D^n}]=0$. By the continuity property of $A$, we have $\hat{\mathbb{E}}[-A_D]=0$. Then it is easy to check that
\begin{displaymath}
Y_0=\hat{\mathbb{E}}[\int_0^D f(s,Y_s,Z_s)ds+S_D I_{\{D<T\}}+\xi I_{\{D=T\}}].
\end{displaymath}
Hence, the result follows.
\end{proof}

\appendix
\renewcommand\thesection{Appendix}
\section{ }
\renewcommand\thesection{A}
In this section, we mainly introduce the extended conditional $G$-expectation and optional stopping theorem under $G$-framework. More details can be found in \cite{HP13}.

Let $(\Omega,L_{ip}(\Omega),\hat{\mathbb{E}}[\cdot])$ be the $G$-expectation space and $\mathcal{P}$ be a weakly compact set that represents $\hat{\mathbb{E}}$. We set
\begin{align*} &L^0(\Omega):=\{X:\Omega\rightarrow [-\infty,\infty] \textrm{ and } X \textrm{ is } \mathcal{B}(\Omega)\textrm{-measurable}\},\\
&\mathcal{L}(\Omega):=\{X\in L^0(\Omega):E_P[X] \textrm{ exists for each }P\in \mathcal{P}\}.
\end{align*}

We extend $G$-expectation $\hat{\mathbb{E}}$ to $\mathcal{L}(\Omega)$ and still denote it by $\hat{\mathbb{E}}$. For each $X\in\mathcal{L}(\Omega)$, we define
\[\hat{\mathbb{E}}[X]=\sup_{P\in\mathcal{P}}E_P[X].\]

Then we give some notations
\begin{align*}
\mathbb{L}^p(\Omega)&:=\{X\in L^0(\Omega):\hat{\mathbb{E}}[|X|^p]<\infty\} \textrm{ for } p\geq1,\\
L_G^{1^*}(\Omega)&:=\{X\in\mathbb{L}^1(\Omega):\exists X_n\in L_G^1(\Omega) \textrm{ such that } X_n\downarrow X, q.s.\},\\
L_G^{1^*_*}(\Omega)&:=\{X\in\mathbb{L}^1(\Omega):\exists X_n\in L_G^{1^*}(\Omega) \textrm{ such that } X_n\uparrow X, q.s.\},\\
\bar{L}_G^{1^*_*}(\Omega)&:=\{X\in\mathbb{L}^1(\Omega):\exists X_n\in L_G^{1^*_*}(\Omega) \textrm{ such that } \hat{\mathbb{E}}[|X_n-X|\rightarrow 0]\}.
\end{align*}

Set $\Omega_t=\{\omega_{\cdot\wedge t}:\omega\in\Omega\}$ for $t>0$. Similarly, we can define $L^0(\Omega_t)$, $\mathcal{L}(\Omega_t)$, $\mathbb{L}^p(\Omega_t)$, $L_G^{1^*}(\Omega_t)$, $L_G^{1^*_*}(\Omega_t)$ and $\bar{L}_G^{1^*_*}(\Omega_t)$ respectively. Then we can extend the conditional $G$-expectation to space $\bar{L}_G^{1^*_*}(\Omega)$.
\begin{proposition}(\cite{HP13})
For each $X\in \bar{L}_G^{1^*_*}(\Omega)$, we have, for each $P\in\mathcal{P}$,
\[\hat{\mathbb{E}}_t[X]={\esssup_{Q\in\mathcal{P}(t,P)}}^P E_Q[X|\mathcal{F}_t],\ \   P\textrm{-a.s.},\]
where $\mathcal{P}(t,P)=\{Q\in P: E_Q[X]=E_P[X], \forall X\in L_{ip}(\Omega_t)\}$.
\end{proposition}

We now give the definition of stopping times under $G$-expectation framework.
\begin{definition}
A random time $\tau :\Omega\rightarrow [0,\infty)$ is called a $*$-stopping time if $I_{\{\tau\geq t\}}\in L_G^{1^*}(\Omega_t)$ for each $t\geq 0$.
\end{definition}

\begin{definition}
For a given $*$-stopping time $\tau$ and $\xi\in \bar{L}_G^{1^*_*}(\Omega)$, we define $\hat{\mathbb{E}}_\tau[\xi]:=M_\tau$, where $M_t=\hat{\mathbb{E}}_t[\xi]$ for $t\geq0$.
\end{definition}

We then give an example of $*$-stopping time.
\begin{example}\label{th13}
Let $(X_t)_{t\in[0,T]}$ be an d-dimensional right continuous process such that $X_t\in L_G^1(\Omega_t)$ for $t\in[0,T]$. For each fixed closed set $F\in \mathbb{R}^d$, we define
\[\tau=\inf\{t\geq0: X_t\notin F\}\wedge T.\]
Then $\tau$ is a $*$-stopping time.
\end{example}

Now we introduce the following optional stopping theorem under $G$-framework.
\begin{theorem}\label{th14}(\cite{HP13})
Suppose that $G$ is non-degenerate. Let $M_t=\hat{\mathbb{E}}_t[\xi]$ for $t\leq T$, $\xi\in L_G^p(\Omega_T)$ with $p>1$ and let $\sigma$, $\tau$ be two $*$-stopping times with $\sigma\leq \tau\leq T$. Then $M_\tau$, $M_\sigma\in \bar{L}_G^{1^*_*}(\Omega)$ and
\[M_\sigma=\hat{\mathbb{E}}_\sigma[M_\tau], \ \ q.s..\]
\end{theorem}

\appendix
\renewcommand\thesection{Acknowledgements}
\section{ }
The authors are grateful to Yongsheng Song for his help and many useful discussions.

\end{document}